\documentclass[12pt]{article}
\usepackage{amsmath, amsthm, amsfonts, amssymb, mathtools}
\usepackage[utf8]{inputenc}
\usepackage[all,cmtip]{xy}
\usepackage{xcolor}
\usepackage{graphicx}
\graphicspath{ {./s/} }
\usepackage{tikz}
\usetikzlibrary{arrows.meta}
\usepackage{float}

\advance \topmargin by -\headheight
\advance \topmargin by -\headsep
\evensidemargin \oddsidemargin
\newtheorem{thm}{Theorem}[section]
\newtheorem{lemma}[thm]{Lemma}
\newtheorem{prop}[thm]{Proposition}
\newtheorem{cor}[thm]{Corollary}

\theoremstyle{definition}
\newtheorem{defn}[thm]{Definition}

\newtheorem{hypoth}{Hypothesis}%[section]

\theoremstyle{remark}
\newtheorem{remark}[thm]{Remark}

\newtheorem{rem}[thm]{Remark}

\numberwithin{equation}{section}

\newcommand{\mbo}{\mathbb{O}}

\DeclareMathOperator{\ev}{ev}
\DeclareMathOperator{\coev}{coev}

\DeclareMathOperator{\Lie}{Lie}
\newcommand{\LLL}{\mathbb{L}}
\DeclareMathOperator{\ad}{ad}

\DeclareMathOperator{\ann}{Ann}

\DeclareMathOperator{\en}{End}

\DeclareMathOperator{\Hom}{Hom}

\DeclareMathOperator{\res}{Res}

\DeclareMathOperator{\tr}{Tr}
\DeclareMathOperator{\zhu}{Zhu}

\DeclareMathOperator{\supp}{Supp}
\DeclareMathOperator{\Var}{Var}

\newcommand{\Fch}{F^{\text{ch}}}
\newcommand{\low}{\text{low}}
\newcommand{\KL}{\text{KL}}
\newcommand{\Cl}{\mathcal{C}\ell}
\newcommand{\Qst}{Q^{\text{st}}}
\newcommand{\heit}{\text{ht}}

\newcommand{\C}{\mathbb{C}}
\newcommand{\F}{\mathbb{F}}
\newcommand{\R}{\mathbb{R}}
\newcommand{\Q}{\mathbb{Q}}
\newcommand{\Z}{\mathbb{Z}}

\newcommand{\CC}{\mathcal{C}}

\newcommand{\HH}{\mathcal{H}}

\newcommand{\OO}{\mathcal{O}}

\newcommand{\CF}{\mathcal{F}}

\newcommand{\ma}{\mathfrak{a}}
\newcommand{\g}{\mathfrak{g}}
\newcommand{\h}{\mathfrak{h}}

\newcommand{\n}{\mathfrak{n}}

\newcommand{\al}{\alpha}

\newcommand{\G}{\Gamma}
\newcommand{\D}{\Delta}

\newcommand{\la}{\lambda}
\newcommand{\La}{\Lambda}

\newcommand{\ov}{\overline}

\newcommand{\vac}{{\left|0\right>}}

\newcommand{\vir}{\text{Vir}}

\DeclareMathOperator{\im}{Im}

\DeclareMathOperator{\gr}{gr}

\newcommand{\what}[1]{\widehat{#1}}
\newcommand{\wtil}[1]{\widetilde{#1}}

\begin{document}

\begin{center}
{\LARGE \bf Notes on affine $W$-algebras} \par \bigskip

\renewcommand*{\thefootnote}{\fnsymbol{footnote}}
{\normalsize
	Jethro van Ekeren\footnote{\texttt{jethro@impa.br}}
}

\par \bigskip
%
%
%\textsuperscript{1}{\footnotesize Instituto de Matem\'{a}tica e Estat\'{i}stica (GMA), UFF, Niter\'{o}i RJ, Brazil}
%
%\par
%
{\footnotesize Instituto de Matem\'{a}tica Pura e Aplicada (IMPA)} \\ \footnotesize{Rio de Janeiro, Brazil}

\par \bigskip
\end{center}

\vspace*{10mm}

\noindent
\textbf{Abstract.} These are expanded and revised notes for a minicourse entitled ``Affine $W$-algebras'', which took place as part of the thematic month ``Quantum Symmetries'' at the Centre de Recherches Math\'{e}matiques in Montr\'{e}al, Canada in October 2022. The first few sections consist of rapid introductions to vertex algebras, affine Kac-Moody algebras and their integrable and admissible modules, and the homological BRST procedure for quantum Hamiltonian reduction. The affine $W$-algebras are defined using these ingredients. The remainder of the notes treats the structure and representation theory of the exceptional affine $W$-algebras, with emphasis on modular tensor categories of representations.

\vspace*{10mm}

\section{Introduction}

There are various motivations to study affine $W$-algebras: they appear in the analysis of completely integrable models, in conformal field theory, representation theory in general, and the geometric Langlands program in particular. In keeping with the theme of the program, our goal will be to understand something of the interesting modular tensor categories that arise as categories of modules over \emph{exceptional} affine $W$-algebras.

Our story begins with the affine Kac-Moody algebra $\what{\g}$ associated with a finite dimensional simple Lie algebra $\g$. For each non negative integer $k$ there is a finite set of integrable highest weight $\what\g$-modules $L_k(\lambda)$ of level $k$, whose characters are given by Kac's generalisation of the Weyl formula \cite{Kac74}. The formula exhibits the normalised characters as quotients of theta functions, with the remarkable consequence that they are modular functions. The explicit modular transformations were determined by Kac and Peterson \cite{KP84}, so that in particular
\[
\chi_{\la}(-1/\tau) = \sum_{\la'} S_{\la, \la'} \chi_{\la'}(\tau),
\]
where $S$ is the matrix
\begin{align}\label{eq:Kac.sum}
S_{\la, \la'} = \frac{i^{|\D_+|}}{|P/(k+h^\vee)Q^\vee|^{1/2}} \sum_{w \in W} \varepsilon(w) e^{-\frac{2\pi i}{k+h^\vee} \left(w(\la+\rho), \la'+\rho\right)}.
\end{align}
Here $P$ and $Q^\vee$ are the weight lattice and coroot lattice of $\g$, respectively, $h^\vee$ is the dual Coxeter number and $\rho$ the Weyl vector, and sum runs over the Weyl group $W$ of $\g$. The finite index set for the highest weights $\la, \la'$ here is $P_+^k$, defined in \eqref{eq:P_+k.def} below.

It emerged that the integrable $\what\g$-module $L_k(\g) = L_k(0)$ is an example of a new kind of algebraic structure known as a vertex algebra \cite{B}. In general a vertex algebra is said to be rational if its modules are completely reducible into direct sums of irreducible ones, and if there are finitely many of the latter up to isomorphism. In fact the affine vertex algebra $L_k(\g)$ is rational, and the finitely many irreducible $L_k(\g)$-modules are precisely the $\what\g$-modules $L_k(\la)$ considered above \cite{FZ92}.

Theorems of Zhu \cite{zhu96} and Huang \cite{H.MTC} assert, respectively, that if a vertex algebra is rational (and satisfies a few other more technical hypotheses) then appropriately normalised characters of its irreducible modules converge to modular functions, and its category of modules is a modular tensor category whose $S$-matrix coincides with the modular $S$-matrix of its characters. Here the tensor product of modules is similar in spirit to, and generalises, the Kazhdan-Lusztig tensor product \cite{KL1} \cite{HL1}. In particular the decomposition of tensor products can be determined through the Verlinde formula \cite{Verlinde} \cite{MS}, and we see in particular that the formula \eqref{eq:Kac.sum} is key to understanding the tensor product structure of the category of $L_k(\g)$-modules.

Affine $W$-algebras are defined via BRST reduction of affine vertex algebras, relative to unipotent groups associated with nilpotent elements. This reduction, known as quantised Drinfeld-Sokolov reduction \cite{DS}, was introduced by Feigin and Frenkel for the principal nilpotent orbit \cite{FF} and in general by Kac, Roan and Wakimoto \cite{KRW}. Varying the simple Lie algebra $\g$, the level $k$, and the choice of nonzero nilpotent element $f \in \g$, results in a bewildering multitude of vertex algebras. It turns out that the reduction of $L_k(\g)$ for $k \in \Z_+$ vanishes for all nonzero $f$, but interesting examples are obtained when we pass to the Kac-Wakimoto \emph{admissible} levels $k$ \cite{KW88}. There is now a delicate interplay between the dimension of the nilpotent orbit $G \cdot f$, the size of the denominator $q$ of $k$, and the ``size'' of the affine $W$-algebra; for each $q$ there being a distinguished nilpotent orbit $G \cdot f = \mbo_q$ for which the corresponding reduction of $L_k(\g)$ has a chance of being rational. The reduction is known as an exceptional $W$-algebra, and the end goal of these notes is to present results on rationality and representation theory of exceptional $W$-algebras.

Characters of irreducible admissible $\what\g$-modules $L_k(\lambda)$, and their modular transformation properties have been determined by Kac and Wakimoto \cite{KW90}, and from these formulas together with structural results about the reduction functor due to Arakawa \cite{Ar.rep.I, Ar.rep.II} the $S$-matrices and fusion rules of exceptional $W$-algebras can be computed. In the principal case this was already done in \cite{FKW}, where a beautiful factorisation of fusion rules was discovered (see Section \ref{sec:excep.chars} below).

For the $S$-matrices and fusion rules of other exceptional $W$-algebras we uncover curious variants of the sum \eqref{eq:Kac.sum} \cite{AVE-EMS}. For example if $\g$ is a Lie algebra of type $D$ or $E$, and $\al_*$ the root of $\g$ associated with the trivalent node of its Dynkin diagram, then the following matrix appears as the $S$-matrix of a modular tensor category of representations of an exceptional $W$-algebra:
\[
S_{\nu, \nu'} = \sum_{w(\alpha^*) \in \D_+} \varepsilon(w) \frac{\left<w(\alpha_*), x\right>}{\left<\alpha_*, x\right>} e^{-\frac{2\pi i}{q} (y(\nu), \nu')}.
\]
Here the sum runs over half of the Weyl group, and $x \in \h$ is not orthogonal to $\al_*$, but otherwise arbitrary. Again, see Section \ref{sec:excep.chars} below.

We conclude this introduction with a few indications of references. For the basic theory of vertex algebras we recommend the books \cite{Kac.VA} and \cite{FBZ}, as well as \cite{KRR} (where the construction of the Zhu algebra is described). The introduction to \cite{DSK06} also contains a short and clear account of various equivalent definitions. For Kac-Moody algebras we recommend Kac's book \cite{Kac.IDLA}, and for background on category $\OO$ methods, such as translation functors, we recommend \cite{Humphreys}. For BRST reduction in general there is the classic article of Kostant and Sternberg \cite{KS}, and Drinfeld-Sokolov reduction of vertex algebras is described in the book \cite{FBZ} for principal nilpotent $f$. Going further, Arakawa's lecture notes \cite{A.minicourse} are a good source, and beyond this there are the primary sources: articles of Kac and Wakimoto and collaborators \cite{KW88, KW90, FKW, KRW, KW04, KW09} and Arakawa \cite{A.rat.O, Ar.rep.I, Ar.rep.II, Arakawa.2015}.

\emph{Acknowledgements.} I would like to thank the organisers of ``Quantum Symmetries'' for the opportunity to participate in this stimulating program. For the duration of the program I held a CRM-Simons Scholar position, and I am also very grateful for this support. During the writing I have been supported by Serrapilheira Institute grant Serra -- 2023-0001, FAPERJ grant 201.445/2021 and CNPq grant 310576/2020-2.

\section{Introduction to vertex algebras}\label{sec:lec1}

Apart from fixing notation and recording key results and background material, the aim of this section is to dispel the common notion that vertex algebras are difficult to work with and calculate in. To do this we present the $\lambda$-bracket calculus as soon as possible, and use it to do a number of example calculations.

We begin by considering the ring
\[
H = \C[x_1, x_2, x_3, \ldots]
\]
of polynomials in countably many variables. In fact we are not so interested in the ring structure of $H$, instead we consider $H$ as a vector space, and study the following collection of linear endomorphisms of $H$: For each integer $n > 0$ we set
\[
h_{-n} = x_n, \quad \text{and} \quad  h_n = n \frac{\partial}{\partial x_n},
\]
and we set $h_0 = 0$. It is easy to confirm that
\[
h_m h_n - h_n h_m = m \delta_{m, -n} I_H.
\]
Therefore $H$ is a module over the following Lie algebra $\what\ma$, which is known as the oscillator Lie algebra. As a vector space
\[
\what\ma = \C[t^{\pm 1}] \oplus \C K, 
\]
and the Lie bracket is defined by
\[
[f(t), g(t)] = \res_{t} f'(t) g(t) \, K, \quad \text{$K$ central}.
\]
Here the symbol $\res_t$ is the formal residue, which acts by extracting the coefficient of $t^{-1}$.

The representation of $\what\ma$ in $H$ is given by $t^m \mapsto h_m$ for all $m \in \Z$ and $K \mapsto I_H$. This $\what\ma$-module is very well known, and is often referred to as the Fock module, following its appearance in quantum field theory.

The first basic idea of vertex algebras is to manipulate $\{h_n \mid n \in \Z\}$ and similar collections of endomorphisms, using generating functions. In this case we set
\[
h(z) = \sum_{n \in \Z} h_n z^{-n-1},
\]
noting that at this point the peculiar choice of exponent of $z$ can be thought of as purely conventional. Now we can set about recasting the commutation relations between the $h_n$ in terms of $h(z)$. Indeed
\begin{align*}
[h(z), h(w)]
&= \sum_{m, n \in \Z} [h_m, h_n] z^{-m-1} w^{-n-1} \\
&= \sum_{m, n \in \Z} m \delta_{m, -n} I_H z^{-m-1} w^{-n-1} \\
&= \sum_{m \in \Z} m I_H z^{-m-1} w^{m-1}.
\end{align*}
This last expression is reminiscent of the series expansions
\[
i_{z, w} \frac{1}{(z-w)^2} = \sum_{m =0}^\infty m z^{-m-1} w^{m-1} \quad \text{and} \quad i_{w, z} \frac{1}{(z-w)^2} = -\sum_{m =0}^\infty m z^{m-1} w^{-m-1}
\]
convergent, respectively, in the domains $|w| < |z|$ and $|z| < |w|$. Here $i_{z, w}$ and $i_{w, z}$ denote the operations of expansion in positive powers of $w$, and of $z$, respectively. These operations can alternatively be characterised as the unique embeddings of fields
\[
\C((z))((w)) \leftarrow \C((z, w)) \rightarrow \C((w))((z))
\]
compatible with the inclusions of $\C[z^{\pm 1}, w^{\pm 1}]$.

To put the calculation of $[h(z), h(w)]$ into a more general context, it is useful to introduce some notation. The formal delta function is by definition the series
\begin{align}\label{eq:delta.def}
\delta(z, w) = \sum_{m \in \Z} z^{-m-1} w^{m} = i_{z, w} \frac{1}{z-w} - i_{w, z} \frac{1}{z-w}.
\end{align}
We may then write
\begin{align}\label{eq:h.bracket}
[h(z), h(w)]
= \partial_{w} \delta(z, w) I_H.
\end{align}

\begin{rem}
There are good reasons for applying the name ``delta function'' to \eqref{eq:delta.def}. By way of justification we do no more than recall the identity
\[
\delta(x) = \frac{1}{2\pi i} \lim_{\varepsilon \rightarrow 0^+} \left( \frac{1}{x-i\varepsilon} - \frac{1}{x+i\varepsilon} \right)
\]
of distributions on $\R$.
\end{rem}

Generalising the simple observation that
\[
(z-w)\delta(z, w) = 0,
\]
we have
\[
(z-w)^j \partial_w^{(n)} \delta(z, w)
= \begin{dcases}
\partial_w^{(n-j)} \delta(z, w) & j \leq n \\
0 & j > n,
\end{dcases}
\]
(where $\partial^{(j)}$ denotes $\partial^j / j!$). In particular \eqref{eq:h.bracket} implies
\[
(z-w)^2 [h(z), h(w)] = 0.
\]
This is abstracted into the following fundamental notion:
\begin{defn}
A formal series in two variables
\[
D(z, w) = \sum_{m, n \in \Z} D_{m, n} z^{-m-1} w^{-n-1}
\]
is said to be local if there exists a positive integer $N$ such that
\[
(z-w)^N D = 0.
\]
\end{defn}
The following proposition shows that the existence of expansions like \eqref{eq:h.bracket} is a general property of locality.
\begin{prop}\label{prop:OPE}
The formal series $D(z, w)$ satisfies $(z-w)^N D(z, w) = 0$ if and only if it is expressible in the form
\[
D(z, w) = \sum_{j=0}^{N-1} c_j(w) \partial_w^{(j)} \delta(z, w),
\]
for some collection of formal series $c_j(w)$ in one variable.
\end{prop}
There is actually a simple formula for the coefficients $c_j(w)$ in Proposition \ref{prop:OPE}. If $D(z, w)$ is a local series then
\[
c_j(w) = \res_z (z-w)^j D(z, w).
\]

We now say that a pair of series $a(w)$ and $b(w)$ forms a local pair if $[a(z), b(w)]$ is local. In this sense our example of $h(w)$ forms a local pair with itself. In general, a local pair $(a(w), b(w))$ of series yields a finite list of new series $c_j(w)$ via the expansion
\[
[a(z), b(w)] = \sum_{j=0}^{N-1} c_j(w) \partial^{(j)} \delta(z, w).
\]
The expansion on the right hand side is known as an ``operator product expansion'', or OPE.

We would like to adopt a perspective in which the fundamental objects are the series themselves, and suppress from our notation the explicit dependence on variables $z$ and $w$. We thus rewrite the OPE above, purely formally, as
\[
[a_\lambda b] = \sum_{j=0}^{N-1} c_j \lambda^{(j)}.
\]
%This is called the OPE of $a(w)$ and $b(w)$. More precisely, we write $a(w)_{(j)}b(w)$ for the new quantum field $c_j(w)$, and we put them all together into the following ``$\lambda$-bracket'':
For instance if we write $h = h(w)$ and $\mathbf{1} = I_H$ then we have
\[
[h_\lambda h] = \lambda \mathbf{1}.
\]

To get to the definition of vertex algebra, one further ingredient is required. 
\begin{defn}
Let $V$ be a vector space. We say that the formal series
\[
a(w) = \sum_{n \in \Z} a_{(n)} w^{-n-1}, \qquad a_{(n)} \in \en(V),
\]
is a quantum field on $V$ if, for every $v \in V$, there exists $N$ such that $a_{(n)}v = 0$ for all $n > N$. In other words if $a(z)v \in V((z))$.
\end{defn}
It is easy to see that $h(w)$ is a quantum field on $H$.

The significance of the quantum field property is that it permits the definition of a formal analogue of the normal ordering procedure of quantum field theory. Specifically, for $a(w)$ and $b(w)$ quantum fields on a vector space $V$, their normally ordered product is by definition
\begin{align*}
:a(w) b(w): = a(w)_+ b(w) + b(w) a(w)_-
\end{align*}
where
\[
a(w)_+ = \sum_{n < 0} a_{(n)} w^{-n-1} \quad \text{and} \quad a(w)_- = \sum_{n \geq 0} a_{(n)} w^{-n-1}.
\]
\begin{prop}
If $a(w)$ and $b(w)$ are two quantum fields on a vector space $V$, then the normally ordered product $:a(w)b(w):$ is a well defined and is another quantum field on $V$.
\end{prop}
Happily the quantum field property is also preserved upon passage to OPE coefficients.
\begin{prop}
If $a(w)$ and $b(w)$ are two quantum fields on a vector space $V$ which form a local pair, then the series $a(w)_{(j)}b(w)$, for $j \geq 0$, is also a quantum field on $V$.
\end{prop}

In fact the two constructions are more similar than they look. It is straightforward to verify 
\begin{align*}
:a(w) b(w): = \res_z \left( a(z) b(w) i_{z, w}(z-w)^{-1} - b(w) a(z) i_{w, z} (z-w)^{-1} \right),
\end{align*}
so that it makes sense to write $:a(w)b(w): = a(w)_{(-1)}b(w)$. In fact the product
\begin{align*}
a(w)_{(n)}b(w) = \res_z \left( a(z) b(w) i_{z, w}(z-w)^{n} - b(w) a(z) i_{w, z} (z-w)^{n} \right)
\end{align*}
makes sense for any pair of quantum fields $a(w)$ and $b(w)$, and is again a quantum field. On the other hand these $n^{\text{th}}$ products, for $n \leq -2$, can easily be recovered from the normally ordered product together with the natural operation of differentiation of a quantum fields:
\begin{align}\label{eq:deriv.nth}
(\partial_w a(w))_{(n)} b(w) = -n a(w)_{(n-1)} b(w).
\end{align}

The relation \eqref{eq:deriv.nth} implies the following relation satisfied by the $\la$-bracket:
\begin{align}\label{eq:la.id.1}
[(\partial a)_{\la}b] = -\la [a_\la b].
\end{align}
Other relations such as
\begin{align}\label{eq:la.id.2}
[a_\la (\partial b)] = (\partial + \lambda) [a_\la b] \quad \text{and} \quad :ab: = :(\partial a)b: + :a(\partial b):
\end{align}
are also easily established.

We would also like formulas allowing us to compute $\lambda$-brackets of normally ordered products of quantum fields, and indeed to manipulate arbitrary combinations of $\lambda$-brackets, normally ordered products and derivatives. Vertex algebras provide a context in which we can do this.
\begin{defn}\label{def:VA.1}
A vertex algebra consists of a vector space $V$, an operator $T : V \rightarrow V$, a nonzero vector $\vac$ such that $T\vac = 0$, and a set $\CF$ of quantum fields on $V$, such that
\begin{itemize}
\item For all $a(w) = \sum_{n \in \Z} a_{(n)} w^{-n-1} \in \CF$, we have
\[
[T, a(w)] = \partial_w a(w),
\]

\item All the quantum fields in $\CF$ are mutually local,

\item $V$ is generated from $\vac$ by $\CF$, in the sense that the terms
\[
a^1_{(n_1)} \cdots a^k_{(n_k)} \vac
\]
together span $V$, as $a^1(w), \ldots, a^k(w)$ run over $\CF$ and $n_1, \ldots, n_k$ run over $\Z$.
\end{itemize}
\end{defn}
The notions of quantum field and mutual locality are incorporated into the definition of vertex algebra. The axiom about generation of $V$ by the coefficients of elements of $\CF$ is inessential in the sense that it can be achieved by replacing $V$ by the generated subspace. The two further ingredients in the definition are the vacuum vector $\vac$, and the translation operator $T$.

It is now possible to set up a ``calculus of quantum fields'' as outlined above, for the elements of $\CF$ and their iterated $\lambda$-brackets and normally ordered products. Alongside the formulas \eqref{eq:la.id.1} and \eqref{eq:la.id.2}, we can also prove
\begin{align*}
[b_\lambda a] = -[a_{-\lambda-\partial} b]
\end{align*}
and
\begin{align*}
[a_\la :bc:] = :[a_\la b]c: + :b[a_\la c]: + \int_0^\la [[a_\la b]_\mu c] \, d\mu.
\end{align*}

Let us return to our original example $H$, check that it is indeed a vertex algebra, and do some computations. We equip the vector space $H$ with $\vac = 1$ and $\CF = \{h(w)\}$. The operator $T : H \rightarrow H$ is determined uniquely by the conditions $T\vac=0$ and $[T, h(w)] = \partial_w h(w)$, and is well defined since its action can be written down explicitly relative to the basis of $H$ consisting of monomials. It is an easy exercise to check that
\[
T = \sum_{j \in \Z_{\geq 1}} h_{-1-j} h_{j}.
\]
So we have a vertex algebra $H$, which we shall refer to as the Heisenberg vertex algebra.

As an example computation in $H$, let us consider the quantum field
\[
L = \frac{1}{2} :hh:
\]
and try to compute $[L_\la L]$. Using the rules given above, we confirm that
\[
[L_\la h] = \partial h + \la h,
\]
and proceed from there to compute
\[
[L_\la L] = \partial L + 2\la L + \frac{\la^3}{12}.
\]

Upon decoding the $\lambda$-bracket notation, we find
\begin{align}\label{eq:vir.field.bracket}
[L(z), L(w)] = \partial_w L(w) \delta(z, w)  + 2 L(w) \partial_w \delta(z, w) + \frac{1}{12} \partial_w^{3} \delta(z, w),
\end{align}
and hence, setting $L(z) = \sum_{n} L_n z^{-n-2}$, we obtain the relations
\begin{align}\label{eq:vir.rel}
[L_m, L_n] = (m-n) L_{m+n} + \delta_{m, -n} \frac{m^3-m}{12} I_H
\end{align}
between the endomorphisms $L_n \in \en(H)$. In fact we have just constructed a representation in $H$ of the Virasoro Lie algebra, which is by definition
\[
\vir = \bigoplus_{n \in \Z} \C L_n \oplus \C C
\]
with Lie brackets given by $C$ central, and \eqref{eq:vir.rel} with $C$ in place of $I_H$. In general a $\vir$-module $M$ in which $C \mapsto c I_M$ is said to be of central charge $c$, so the $\vir$-module $H$ is of central charge $1$.

We leave it as an exercise to check that the quantum field $B = L + \beta \partial h$ satisfies
\[
[B_\la B] = \partial B + 2\la B + \frac{\la^3}{12}(1-12\beta^2),
\]
so its coefficients endow $H$ with a representation of $\vir$ of central charge $1-12\beta^2$.

If we examine the definition of the field $L(w)$ we see that its coefficients $L_n$ are given by the formulas
\[
L_n = \frac{1}{2} \sum_{j \in \Z} h_{n-j} h_{j}, \quad n \neq 0, \quad \text{and} \quad L_0 = \frac{1}{2} h_0^2 + \sum_{j \in \Z_{\geq 1}} h_{-j} h_j.
\]
Thus the conversion of $H$ from a $\what\ma$-module into a $\vir$-module did not come from anything like a homomorphism between these Lie algebras, we had to step outside the world of Lie algebras to achieve it.

Let us remark at this point, that the coefficient $L_{-1} : H \rightarrow H$ of our Virasoro quantum field $H$ happens to coincide with the translation operator $T : H \rightarrow H$.

Definition \ref{def:VA.1} is a little peculiar for what purports to be the definition of a class of algebras. But it can be transformed into an equivalent but rather different looking form. We describe the reformulation briefly in the following paragraphs.

Let $(V, \vac, T, \CF)$ be a vertex algebra. A result known as Dong's lemma guarantees that locality is preserved by $n^{\text{th}}$ products. We may thus form the envelope $\ov\CF$ of $\CF$ by repeatedly adding all $n^{\text{th}}$ products of existing quantum fields, and pass to the ``completed'' vertex algebra $(V, \vac, T, \ov\CF)$.

The condition $[T, a(w)] = \partial_w a(w)$ on a quantum field $a(w)$ guarantees that $a(w)\vac \in V[[w]]$. In particular there is a well defined linear map
\[
a(w) \mapsto a(w)\vac|_{w=0}, \qquad \ov{\CF} \rightarrow V.
\]
A key result in the basic theory of vertex algebras is that this map is a bijection. We thus denote by $a \mapsto Y(a, w)$ the inverse of this map, then it turns out that
\[
Y(\vac, w) = I_V, \quad \text{and} \quad Y(a, w)_{(n)}Y(b, w) = Y(a_{(n)}b, w).
\]
Expanding out the latter identity in powers of $w$ yields the following Borcherds identity, valid in any vertex algebra $V$, for all  $a, b, c \in V$ and all $m, k, n \in \Z$;
\begin{align}\label{eq:bor.def}
\begin{split}
& \sum_{j \in \Z_+} \binom{m}{j} (a_{(n+j)}b)_{(m+k-j)} c \\
& \,\,\,\,\,\,\,\, = \sum_{j \in \Z_+} (-1)^j \binom{n}{j} \left( a_{(m+n-j)} b_{(k+j)} c - (-1)^n b_{(n+k-j)} a_{(m+j)} c \right).
\end{split}
\end{align}

Finally a vertex algebra can equivalently be defined in terms of the Borcherds identity, i.e.,
\begin{defn}
A vertex algebra is a vector space $V$, with a nonzero vector $\vac \in V$ and a collection of bilinear operations $V \times V \rightarrow V$, indexed by $\Z$, denoted $a_{(n)}b$, such that
\begin{itemize}
\item for all $a, b \in V$ there exists $N$ such that $a_{(n)}b = 0$ for all $n \geq N$,

\item writing $Y(a, w) = \sum_{n \in \Z} a_{(n)} w^{-n-1}$, we have $Y(a, w)\vac \in V[[w]]$ for all $a \in V$,

\item $Y(\vac, w) = I_V$ and $Y(a, w)\vac|_{w=0} = a$, and

\item for all $a, b, c \in V$ and $n, m, k \in \Z$ the identity \eqref{eq:bor.def} holds.
\end{itemize}
\end{defn}

As for associative algebras, commutative algebras, Lie algebras, etc., there is an obvious extension of the notion of vertex algebra to that of vertex superalgebra. Vertex superalgebra play an important role in the theory, appearing for example as vertex algebra analogues of Clifford algebras.

We now present a simple example of a vertex superalgebra: the charged free fermions. Let $N$ be a finite dimensional vector space, and consider the sum $S = N \oplus N^*$ together with the canonical symmetric bilinear form given by pairing
\[
\left<\varphi, n\right> = \left<n, \varphi\right> = \varphi(n), \qquad \text{for $n \in N$, $\varphi \in N^*$}.
\]
We convert $S$ into a purely odd vector superspace, so that $\left<\cdot, \cdot\right>$ is now a super-skewsymmetric pairing. Consider fields $a = a(z)$, one for each $a \in S$, with $\lambda$-bracket defined to be
\[
[a_\lambda b] = \left<a, b\right>.
\]
More precisely we set
\[
\Fch(N) = U(\what{S}) \otimes_{U({\what{S}}_+)} \C \vac
\]
where
\[
\widehat{S} = t^{1/2}S[t^{\pm 1}] \oplus \C K, \qquad [a t^m, b t^n] = \delta_{m, -n} \left<a, b\right> K, \qquad \text{$K$ central},
\]
and the action of $\widehat{S}_+ = t^{1/2}S[t] \oplus \C K$ on $\C \vac$ is given by $at^m = 0$ for all $a \in S$ and $K=1$. It can be confirmed that $\Fch(N)$ is a vertex superalgebra.

Now suppose that $\sigma$ is a finite dimensional representation of a Lie algebra $\g$ in our vector space $N$, let $\{e_1, \ldots, e_r\}$ be a basis of $N$ and $\{\varphi_1, \ldots, \varphi_r\}$ the dual basis of $N^*$. We form the free fermion vertex algebra $\Fch(N)$ as above, and now for each $x \in \g$ we define a quantum field
\[
F^x = \sum_{i=1}^r :(\sigma(x) e_i) \, \varphi_i:
\]
on $\Fch(N)$. It can be checked by explicit calculation with the $\la$-brackets that
\[
[F^{a}_\lambda F^{b}] = F^{[a, b]} + \lambda \, (a, b)_N,
\]
where
\[
(a, b)_N = \tr_N \sigma(a) \sigma(b)
\]
is the trace form. % {\color{red}Sign and normalisation here seem to be OK; see KW-quantum-reduction paper, formula (2.5).}
If we take the collection of fields $F^{x}(z)$ for $x \in \g$, and the subspace of $\Fch(N)$ generated from $\vac$ by these fields, then we obtain a vertex algebra, canonically associated to $\g$ and its representation in $N$. In Section \ref{sec:lec2} we will examine a general construction of an ``affine'' vertex algebra associated with a Lie algebra $\g$ and an invariant bilinear form on $\g$.

We now very briefly discuss cohomology of vertex superalgebras. Let $V$ be a vertex superalgebra and let $Q \in V$. It follows from the Borcherds identity that the operator
\[
Q_{(0)} : V \rightarrow V
\]
satisfies
\[
[Q_{(0)}, a_{(n)}] = (Q_{(0)}a)_{(n)}, \qquad \text{for all $a \in V$, $n \in \Z$}.
\]

Now, there are natural notions of vertex subalgebra and vertex ideal (and quotient vertex algebra). We observe that $\ker(Q_{(0)}) \subset V$ is a vertex subalgebra. Now suppose that $Q$ is odd and satisfies $Q_{(0)}Q = 0$, then a short computation reveals that $\im(Q_{(0)}) \subset \ker(Q_{(0)})$ is an ideal. An odd element $Q$ as above thus gives a notion of homology algebra $H(V, Q) = \ker(Q_{(0)}) / \im(Q_{(0)})$.

Later we will become familiar with the notion of \emph{conformal vertex algebra}, for now we just note that if $V$ is a conformal vertex algebra, then the operator $Q_{(0)} : V \rightarrow V$ will preserve conformal weight if $\Delta(Q) = 1$.

\section{Affine vertex algebras and the Zhu algebra}\label{sec:lec2}

Let $\g$ be a finite dimensional Lie algebra over $\C$, and $\phi : \g \times \g \rightarrow \C$ an invariant bilinear form. The associated affine Lie algebra is, by definition,
\[
\what\g = \g[t, t^{-1}] \oplus \C K
\]
with Lie bracket
\begin{align}\label{eq:aff.lie.bracket}
[a f, b g] = [a, b] fg + \phi(a, b) \res_t f'(t) g(t) K, \qquad \text{$K$ central}.
\end{align}
The oscillator algebra $\what\ma$ is recovered as the special case $\g = \C$ with the bilinear form $\phi(1, 1) = 1$. The Fock space $H$ of $\what\ma$ is generalised by the construction of the vacuum $\what\g$-module as follows. Set
\[
V^\phi(\g) = U(\what\g) \otimes_{U(\g[t]+\C K)} \C\vac.
\]
Here $\C \vac$ is a one dimensional vector space furnished with an action of the subalgebra $\g[t]+\C K \subset \what\g$ by $K \mapsto 1$ and $\g[t] \mapsto 0$. %The scalar $k$ is referred to as the level of the vacuum representation.
The vacuum module $V^\phi(\g)$ is now equipped with quantum fields $x(z) = \sum_{n \in \Z} x t^n \, z^{-n-1}$, one for each $x \in \g$, and these generate a vertex algebra structure.

For $x \in \g$ we write $x_n \in \en(V^\phi(\g))$ as a shorthand for the action of $xt^n$, and we write $x_{-1}\vac \in V^\phi(\g)$ simply as $x$. Comparing these conventions with those of Section \ref{sec:lec1} we see that $x_{(n)} = x_n$ for all $n \in \Z$ for $x \in \g$. See also the remarks following Definition \ref{defn:conf} below.

The construction of fields $F^{x}(z)$ in the free fermions vertex superalgebra $\Fch(N)$ of a $\g$-module $N$, may now be interpreted as the construction of a homomorphism of vertex algebras
\[
V^\phi(\g) \rightarrow \Fch(N),
\]
where $\phi$ is the trace form of $N$.

The construction of $V^\phi(\g)$ is most commonly carried out in the following special case. Let $\g$ be a finite dimensional simple Lie algebra, $(\cdot, \cdot) : \g \times \g \rightarrow \C$ the normalised bilinear form (so that $(\theta, \theta) = 2$ for $\theta$ the highest root of $\g$), and $\phi(x, y) = k \cdot (x, y)$. We then write $V^k(\g)$ without ambiguity.

\begin{figure}[H]
    \centering
\begin{tikzpicture}
  % Grid of dots
  \def\s{2};
  \node[circle,fill,inner sep=2pt, label={$1$}] at (0,0) {};
  \node[circle,fill,inner sep=2pt, label={$f_{-1}$}] at (-\s,-\s) {};
  \node[circle,fill,inner sep=2pt, label={$h_{-1}$}] at (0,-\s) {};
  \node[circle,fill,inner sep=2pt, label={$e_{-1}$}] at (\s,-\s) {};
  \node[circle,fill,inner sep=2pt, label=below:{$f_{-1}^2$}] at (-2*\s,-2*\s) {};
  \node[circle,fill,inner sep=2pt, label={$f_{-1}h_{-1}$}, label=below:{$f_{-2}$}] at (-\s,-2*\s) {};
  \node[circle,fill,inner sep=2pt, label={$h_{-2}$}, label=below:{$f_{-1}e_{-1}\,\,\, h_{-1}^2$}] at (0,-2*\s) {};
  \node[circle,fill,inner sep=2pt, label={$h_{-1}e_{-1}$}, label=below:{$e_{-2}$}] at (\s,-2*\s) {};
  \node[circle,fill,inner sep=2pt, label=below:{$e_{-1}^2$}] at (2*\s,{-2*\s}) {};

  % Dotted axes
  \draw[dotted] (0,0) -- (-6,-6);
  \draw[dotted] (-2,-2) -- (2,-2);
  \draw[dotted] (-4,-4) -- (4,-4);
  \draw[dotted] (0,0) -- (6,-6);
  
\end{tikzpicture}
    \caption{Sketch of the vertex algebra $V^k(\mathfrak{sl_2})$}
    \label{fig:Vkg}
\end{figure}
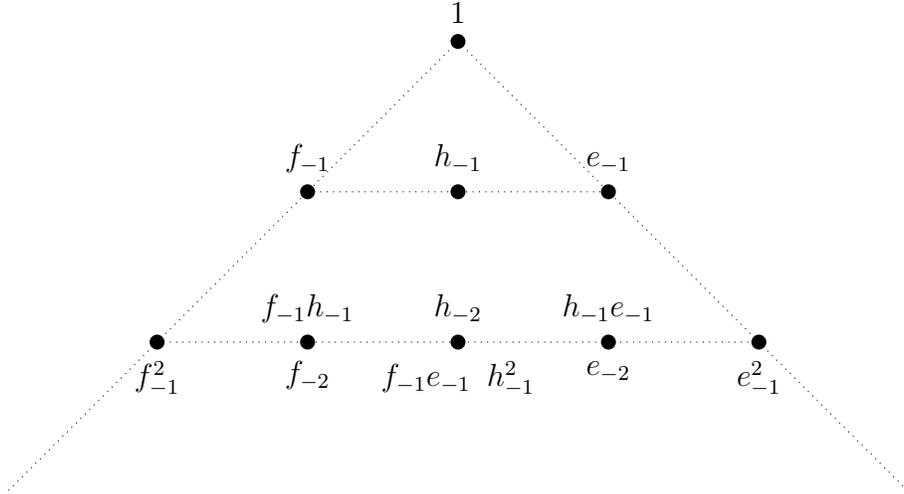

Let $\{a^1, \ldots, a^m\}$ be a basis of $\g$. By the PBW theorem, a basis of $V^k(\g)$ is provided by the monomials
\begin{align}\label{eq:Vkg.monom}
a^{i_1}_{(-n_1)} a^{i_2}_{(-n_2)} \cdots a^{i_s}_{(-n_s)} \vac
\end{align}
where $n_1 \geq n_2 \geq \ldots \geq n_s \geq 1$ and $1 \leq i_k \leq m$ for each $k$, and if $n_k = n_{k+1}$ then $i_k \geq i_{k+1}$. A visual impression of $V^k(\mathfrak{sl}_2)$ appears in Figure \ref{fig:Vkg}.

As usual we let $\mathbf{1}$ denote the quantum field $I_{V}$ on $V = V^k(\g)$. To verify mutual locality of the fields $a(w)$ introduced already, it suffices to convert the commutation relation \eqref{eq:aff.lie.bracket} into a relation between fields $a(z)$ and $b(w)$, which is then seen to be equivalent to the following $\la$-bracket relation
\[
[a_\lambda b] = [a, b] + k (a, b) \lambda \mathbf{1}.
\]

We continue to assume $\g$ is simple. Let $h^\vee$ denote the dual Coxeter number of $\g$, and let $\{a_i\}$ and $\{b_i\}$ be bases of $\g$ dual with respect to $(,)$. If $k \neq -h^\vee$ then the Sugawara construction yields a vector
\[
L = \frac{1}{2(k+h^\vee)} \sum_i :a_i b_i:
\]
which satisfies Virasoro algebra relations, i.e., 
\[
[L_\lambda L] = \partial L + 2\lambda L + \frac{\lambda^3}{12} c \,\, \mathbf{1},
\]
where in this case the central charge is given by
\[
c = \frac{k \dim(\g)}{k+h^\vee}.
\]

Virasoro vectors in vertex algebras turn out to be important in general for the development of the geometric side of the theory.
\begin{defn}\label{defn:conf}
We say the vertex algebra $V$ is \emph{conformal} if it contains a vector $L$, the modes of whose associated quantum field
\[
Y(L, w) = L(w) = \sum_{n \in \Z} L_n w^{-n-2}
\]
satisfy the relations
\[
[L_m, L_n] = (m-n) L_{m+n} + \delta_{m, -n} \frac{m^3-m}{12} c \,\, I_V,
\]
such that $T = L_{-1}$, and such that $V$ is decomposed as $V = \bigoplus_{n \in \Z_+} V_n$ into finite dimensional eigenspaces of $L_0$. The eigenvalue $n$ of an element $a$ of $V_n$ is called the conformal weight of $a$ and denoted $\D(a)$.
\end{defn}
As a straightforward consequence of the axioms of a vertex algebra, we find that
\[
a_{(n)} : V_\alpha \rightarrow V_{\alpha + (\Delta-n-1)}
\]
for $a \in V_\D$. We thus write $a_n = a_{(n + \Delta -1)}$, so that
\[
a_n : V_\alpha \rightarrow V_{\alpha-n}.
\]
The vertex algebra $V^k(\g)$, with its Sugawara vector $L$, is conformal, and the conformal weight grading is the vertical grading implicit in Figure \ref{fig:Vkg} above. More precisely the conformal weight of the vector \eqref{eq:Vkg.monom} is $n_1 + \cdots + n_s$.

The Borcherds identity can be rewritten in terms of the conformal grading on modes, and becomes
\begin{align}\label{eq:bor}
\begin{split}
& \sum_{j \in \Z_+} \binom{m + \Delta(a)-1}{j} (a_{(n+j)}b)_{n+m+k} c \\
& \,\,\,\,\,\,\,\, = \sum_{j \in \Z_+} (-1)^j \binom{n}{j} \left( a_{m+n-j} b_{k+j} c - (-1)^n b_{n+k-j} a_{m+j} c \right).
\end{split}
\end{align}

Now we discuss modules over vertex algebras. To avoid grappling with the several technically inequivalent definitions of $V$-module that arise in the non conformal case, we suppose from the outset that $V$ is conformal.
\begin{defn}
Let $V$ be a conformal vertex algebra. A (positive energy) $V$-module is a direct sum of graded vector spaces of the form
\[
M = \bigoplus_{\alpha - \alpha_0 \in \Z_{\geq 0}} M_\alpha,
\]
equipped with an action
\[
a \mapsto  Y^M(a, z) = \sum_{n \in \Z} a^M_n z^{-n-\Delta(a)}
\]
satisfying \eqref{eq:bor} for all $a, b \in V$ and $c \in M$.
\end{defn}
Let us think about the restriction of \eqref{eq:bor} to the lowest graded piece $c \in M_{\low} = M_{\alpha_0}$. In this piece  we have $a^M_n = 0$ for all $n > 0$. Set $m=1$, $k=0$ and $n = -1$ in Borcherds identity. It becomes
\begin{align*}
\sum_{j \in \Z_+} \binom{\Delta(a)}{j} (a_{(j-1)}b)^M_{0} c
& = a_{0}^M b_{0}^M c.
\end{align*}
This motivates us to define
\[
a * b = \sum_{j \in \Z_+} \binom{\Delta(a)}{j} a_{(j-1)}b
\]
in $V$. Similarly we discover that, if
\[
a \circ b = \sum_{j \in \Z_+} \binom{\Delta(a)}{j} a_{(j-2)}b
\]
then $(a \circ b)^M_0 = 0$ in the restriction to $M_{\low}$. This motivates us to define the quotient
\[
\zhu(V) = \frac{V}{V \circ V}, \qquad \text{with the induced product $a * b$}.
\]

A number of miracles now occur. The subspace $V \circ V$ is really an ideal with respect to the product $*$. The quotient is an associative algebra with unit $\vac$. Furthermore the assignment $M \mapsto M_{\low}$, while not so strong as an equivalence of categories, is bijective between the set of irreducible positive energy $V$-modules, and the set of irreducible left $\zhu(V)$-modules \cite{zhu96}. For $N$ an irreducible left $\zhu(V)$-module we denote by $\LLL(N)$ the unique irreducible positive energy $V$-module for which $\LLL(N)_{\low} \cong N$.

The Zhu algebra of $V^k(\g)$ is isomorphic to the universal enveloping algebra $U(\g)$ (no matter what the value of $k$). Roughly speaking, this is because
\[
a \circ v = a_{(-2)} v + a_{(-1)}v,
\]
so the spanning monomials \eqref{eq:Vkg.monom} of $V^k(\g)$ all reduce to monomials of the form
\begin{align}\label{eq:zhu.monom}
a^{i_1}_{(-1)} a^{i_2}_{(-1)} \cdots a^{i_s}_{(-1)} \vac,
\end{align}
which are then identified with corresponding monomials
\[
a^{i_1} a^{i_2} \cdots a^{i_s} \in U(\g).
\]
Actually things are not quite so simple: to achieve an isomorphism $\zhu(V^k(\g)) \rightarrow U(\g)$ not just of vector spaces, but of algebras, the monomial \eqref{eq:zhu.monom} must be identified with $a^{i_s} a^{i_{s-1}} \cdots a^{i_1}$. It is a good exercise to verify this, using \eqref{eq:bor} with $m=0$, $k=1$. See \cite{FZ92}.

\begin{defn}
A $\widehat{\g}$-module $M$ is said to be smooth if for all $a \in \g$ and $x \in M$, there exists $N \in \Z_+$ such that $at^n \cdot x = 0$ for all $n \geq N$.
\end{defn}
If we wish to try and build a $V^k(\g)$-module structure on a $\widehat{\g}$-module $M$, by setting $a(z) = \sum_{n \in \Z} at^n z^{-n-1}$, then clearly it is a prerequisite that $M$ be smooth. On the other hand it is sufficient. See \cite{KRR} (also \cite{Li.local.sys}) for the following theorem and its proof.
\begin{thm}
There is an equivalence between the category of $V^k(\g)$-modules and the category of smooth $\widehat{\g}$-modules of level $k$.
\end{thm}
If the level $k$ is irrational then $V^k(\g)$ is simple, and that's about all there is to the story. If, on the other hand, $k$ is rational, then it is possible to pass to a simple quotient vertex algebra
\[
L_k(\g) = V^k(\g) / N_k.
\]
Here $N_k \subset V^k(\g)$ is the maximal proper vertex algebra ideal which, in this case, is the same thing as the maximal proper $\widehat{\g}$-submodule. It is a very interesting question to ask which $V^k(\g)$-modules descend to $L_k(\g)$-modules. If we content ourselves with the irreducible positive energy modules, then it is enough to work with the Zhu algebra, and in general we have
\[
\zhu(L_k(\g)) = \zhu(V^k(\g) / N_k) = U / I_k
\]
where $I_k$ denotes the two-sided ideal in $\zhu(V)$ generated by the image of $N_k$. 

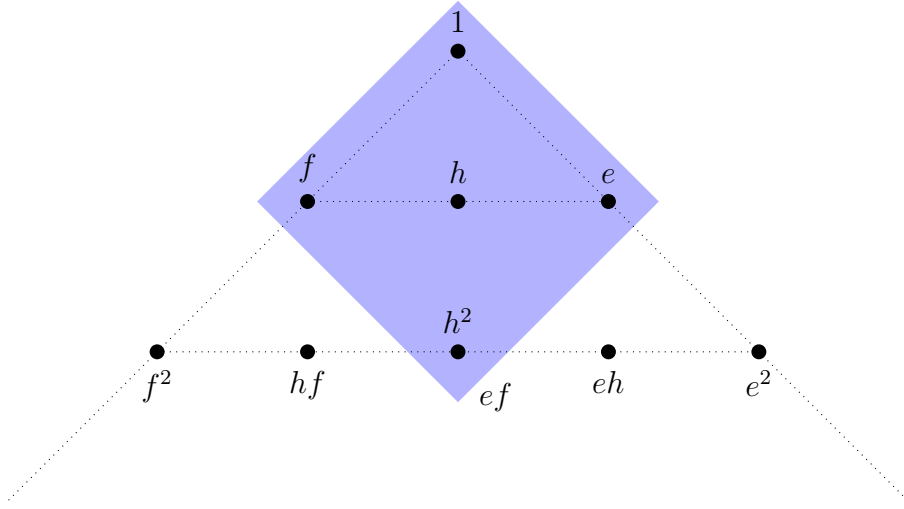
\begin{figure}[H]
    \centering
\begin{tikzpicture}
  \begin{scope}[opacity=0.4]
    \fill[blue!30] (0,0.67) -- (2.67,-2) -- (0,-4.67) -- (-2.67,-2) -- cycle;
  \end{scope}
  % Grid of dots
  \def\s{2};
  \node[circle,fill,inner sep=2pt, label={$1$}] at (0,0) {};
  \node[circle,fill,inner sep=2pt, label={$f$}] at (-\s,-\s) {};
  \node[circle,fill,inner sep=2pt, label={$h$}] at (0,-\s) {};
  \node[circle,fill,inner sep=2pt, label={$e$}] at (\s,-\s) {};
  \node[circle,fill,inner sep=2pt, label=below:{$f^2$}] at (-2*\s,-2*\s) {};
  \node[circle,fill,inner sep=2pt, label=below:{$hf$}] at (-\s,-2*\s) {};
  \node[circle,fill,inner sep=2pt, label={$h^2$}] at (0,-2*\s) {};
  % label=below:{$ef$}
  \node[] at (0.5,-4.6) {$ef$};
  \node[circle,fill,inner sep=2pt, label=below:{$eh$}] at (\s,-2*\s) {};
  \node[circle,fill,inner sep=2pt, label=below:{$e^2$}] at (2*\s,{-2*\s}) {};
  
  % Dotted axes
  \draw[dotted] (0,0) -- (-6,-6);
  \draw[dotted] (-2,-2) -- (2,-2);
  \draw[dotted] (-4,-4) -- (4,-4);
  \draw[dotted] (0,0) -- (6,-6);
  
\end{tikzpicture}
    \caption{Sketch of $U$ with a basis of $U / (e^2)$ highlighted.}
    \label{fig:Ug.k=1}
\end{figure}

For instance in $V^1(\mathfrak{sl}_2)$ the maximal proper ideal is generated by $e_{(-1)}^2\vac$. Writing $U = U(\mathfrak{sl}_2)$, we obtain
\[
\zhu(L_1(\mathfrak{sl}_2)) \cong U / (e^2).
\]
It is now not hard to compute $U/(e^2) \cong \C \oplus \text{Mat}_{2 \times 2}(\C)$ as associative $\C$-algebras. Indeed $e^2 \in I$ implies
\[
[f, e^2] = -(he + eh) = -2eh - 2e \in I,
\]
and after applying $\ad(f)$ again we discover that $2ef - (h^2+h) \in I$, etc. Continuing in this way we eventually conclude that all elements of $U$ can be reduced modulo $I$ to combinations of the $5$ basis elements indicated in Figure \ref{fig:Ug.k=1}. We can also compute the products of these elements explicitly.
%Notice that
%\[
%(ef + fe)^2 = ... = h^2 = ef+fe
%\]
%is an idempotent. It equals the identity on the two-dimensional irreducible $U$-module $L(\varpi_1)$
The algebra $U/(e^2)$ has two irreducible modules; the one dimensional trivial $\mathfrak{sl}_2$-module $L(0)$, and the two dimensional irreducible $\mathfrak{sl}_2$-module $L(\varpi)$ of highest weight $\varpi$. Here $\varpi$ is the fundamental weight dual to the single simple root $\alpha$ of $\mathfrak{sl}_2$.
%\begin{center}
%\begin{tikzpicture}
%  \node[circle,fill,inner sep=2pt] (dot1) at (0,0) {};
%  \node[circle,fill,inner sep=2pt] (dot2) at (2,0) {};
%
%  \draw[->] (dot1) to [bend left=30] node[above] {$e$} (dot2);
%  \draw[->] (dot2) to [bend left=30] node[below] {$f$} (dot1);
%\end{tikzpicture}
%\end{center}
\begin{center}
\begin{tikzpicture}
  \node[circle,fill,inner sep=2pt] (dot1) at (-4,0) {};
  \node[circle,fill,inner sep=2pt] (dot1) at (0,0) {};
  \node[circle,fill,inner sep=2pt] (dot2) at (2,0) {};

  \draw[->] (0.2,0.2) to [bend left=30] node[above] {$e$} (1.8,0.2);
  \draw[->] (1.8,-0.2) to [bend left=30] node[below] {$f$} (0.2,-0.2);
\end{tikzpicture}
\end{center}
There are thus two irreducible positive energy $L_1(\mathfrak{sl}_2)$-modules. Since $L_1(\mathfrak{sl}_2)$ is simple, it is among these, and indeed
\[
\LLL(L(0)) = L_1(\mathfrak{sl}_2).
\]
Continuing along these lines, the following important theorem is established.
\begin{thm}[{\cite[Theorem 3.1.3]{FZ92}}]\label{thm:aff.k.Zplus}
Let $k \in \Z_+$. The complete set of irreducible positive energy $L_k(\g)$-modules is given by $\LLL(L(\lambda))$ as $\lambda$ ranges over
\begin{align}\label{eq:P_+k.def}
P_+^k = \{\lambda \in P_+ \mid \left<\lambda, \theta^\vee \right> \leq k\}.
\end{align}
%The value of $L_0$ (conformal dimension) on $L_k(\lambda)_{low}$ is
%\[
%h_\lambda = \frac{|\lambda+\rho|^2 - |\rho|^2}{2(k+h^\vee)}.
%\]
\end{thm}
The modules $\LLL(L(\la))$ are irreducible $\what\g$-modules, and it is easy to see they are highest weight $\what\g$-modules in fact. We discuss the theory of highest weight modules in greater depth in the next section.

\section{Kac-Moody algebras and admissible weights}

Let $\g$ be a finite dimensional simple Lie algebra of rank $\ell$. Strictly speaking the affine Lie algebra $\widehat{\g}$ is not quite a Kac-Moody algebra, while its extension
\[
\widetilde{\g} = \widehat{\g} \oplus \C d, \qquad [d, at^n] = n at^n, \quad [d, K] = 0,
\]
is. The Cartan subalgebra of $\widetilde{\g}$ is $\h \oplus \C K \oplus \C d$ and, roughly speaking, the role of the extra vector $d$ is to enlarge the dual ${\what\h}^*$ enough for it to contain the root system of $\what\g$. Setting $\delta$ and $\Lambda_0$ as dual vectors to $d$ and $K$, respectively, the standard invariant bilinear form $(,)$ on $\h^*$ extends to $\wtil\h^*$ with
\[
(\delta, \delta) = (\Lambda, \Lambda) = 0, \quad (\delta, \Lambda_0) = 1.
\]
The root system of $\what\g$ (or rather, of $\wtil\g$) is
\[
\what\D = \{n \delta + \alpha \mid n \in \Z, \,\, \alpha \in \D\} \cup \{n \delta \mid n \in \Z\},
\]
with the former and latter sets of roots being referred to as \emph{real} and \emph{imaginary}, respectively. Positive roots of $\wtil\g$ are those of the form $n \delta + \alpha$ where either $n > 0$ or else $n = 0$ and $\alpha \in \D_+$. Simple roots are defined in the usual way, and the set of simple roots turns out to consist of the simple roots $\alpha_1, \ldots, \alpha_\ell$ of $\g$, together with $\alpha_0 = -\theta + \delta$, where $\theta \in \D_+$ denotes the highest root of the finite root system.

The norm of the root $\alpha + n \delta$ is the same as that of $\alpha \in \D$, and the set of real coroots is
\[
\what\D^{\vee, \text{re}} = \{\nu^{-1}(\alpha) / (\alpha, \alpha) \mid \text{$\alpha \in \what\D$ real}\}.
\]

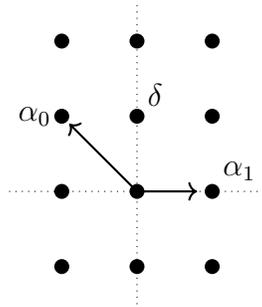
\begin{figure}[h]
    \centering
\begin{tikzpicture}
  % Grid of dots
  \foreach \x in {-1,0,1}
    \foreach \y in {-1,0,1,2}
      \node[circle,fill,inner sep=2pt] at (\x,\y) {};

  % Dotted axes
  \draw[dotted] (-1.7,0) -- (1.7,0);
  \draw[dotted] (0,-1.5) -- (0,2.5);

  % Thick arrows
  \draw[->,thick] (0,0) -- (0.8,0);
  \draw[->,thick] (0,0) -- (-0.9,0.9); 

  % Label (0,1)
  \node[left] at (-1, 1) {$\alpha_0$};
  \node[above right] at (1,0) {$\alpha_1$};
  \node[above right] at (0,1) {$\delta$};
\end{tikzpicture}

    \caption{The root system of $\what{\mathfrak{sl}}_2$}
    \label{fig:sl2hat.rootsystem}
\end{figure}

In Figure \ref{fig:sl2hat.rootsystem} part of the root system of $\what{\mathfrak{sl}}_2$ is displayed. If $\{e, h, f\}$ is the standard basis of $\mathfrak{sl}_2$, then the root spaces associated with $\alpha_1, \alpha_0$ and $\delta$ are spanned by $e$, $f t$ and $ht$, respectively, etc.

We shall generally use the notation $\alpha$ for roots of $\g$ and of $\what\g$, relying on context to resolve ambiguity. On the other hand, having fixed a level $k$, we shall often write $\what\la$ for the affine weight $\what\la = k\Lambda_0 + \la$ corresponding to $\la \in \h^*$. For a $\g$-module $M$ we denote by $V_k(M)$ the corresponding induced $\what\g$-module of level $k$, i.e.,
\[
V_k(M) = U(\what\g) \otimes_{U(\what\g_+)} M,
\]
where $\what\g_+ = \g[t] \oplus \C K$, and its action on $M$ is given by extending that of the subalgebra $\g$ by setting $a t^{>0} = 0$ and $K = k$. If $M$ is a highest weight module of highest weight $\la \in \h^*$ then $V_k(M)$ lies between the corresponding affine Verma module $M(\what\lambda)$ and its irreducible quotient $L(\what\la)$,
\[
M(\what\la) \twoheadrightarrow V_k(M) \twoheadrightarrow L(\what\la).
\]
In particular if $M = M(\la)$ is itself a Verma module, then $V_k(M) \cong M(\what\la)$. At the other extreme, if $M = L(\la)$ is irreducible, then $V_k(M) \cong L(\what\la)$ if the level $k$ is irrational, but in general the surjection can be nontrivial. If $M = L(\la)$ is finite dimensional irreducible, then $V_k(M)$ is known as a Weyl module. We shall denote it by  $V_k(\lambda)$.  Thus for example $V^k(\g) = V_k(0)$ as a $\what\g$-module.  Clearly $V_k(M)_{\low} = M$ in general.

Much as in the finite case, an affine weight $\lambda$ is said to be dominant integral if
\[
\left<\lambda+\widehat\rho, \al\right> \in \Z_{>0} \quad \text{for all $\alpha^\vee \in \widehat\D^{\vee, \text{re}}_+$}.
\]
Here the Weyl vector $\what\rho \in \what\h^*$ is characterised by the property $\what\rho(\alpha_i^\vee) = 1$ for $i=0, 1, \ldots, \ell$. Unlike the finite case, in the affine case this does not characterise $\what\rho$ uniquely, so we make the choice $\what\rho = h^\vee \Lambda_0 + \rho$ for definiteness.

An interesting generalisation of dominant integral in the affine case was identified by Kac and Wakimoto in \cite{KW88}. The integral coroot system of an affine weight $\lambda$ is by definition
\[
\what{\Delta}^{\vee, \text{re}}(\lambda) = \{\al^\vee \in \what\D^{\vee, \text{re}} \mid \left<\lambda+\what\rho, \al^\vee \right> \in \Z\}.
\]
Now $\lambda$ is said to be \emph{admissible} if, in essence, it is dominant integral with respect to its integral coroot system. More precisely, if
\begin{itemize}
\item the rank of $\what{\Delta}^{\vee, \text{re}}(\lambda)$ equals that of $\widehat{\Delta}^\vee$, and

\item $\left<\lambda + \what\rho, \al^\vee \right> \in \Z_{> 0}$ for all $\alpha^\vee \in \what{\Delta}^{\vee, \text{re}}(\lambda) \cap \what{\Delta}_+$.
\end{itemize}

As a simple example, we examine for which $k$ the vacuum weight $k\Lambda_0$ is admissible. The simple coroots in $\what{\D}^{\vee, \text{re}}(k\La_0)$ must be of the form
\[
\Pi^\vee(\lambda) = \{\al^\vee \in \Delta^{\vee, \text{re}}(\lambda) \cap \what{\Delta}_+ \mid \text{$\al$ indecomposable in $ \Delta^{\vee, \text{re}}(\lambda) \cap \what{\Delta}_+$} \}
\]
must be
\[
\Pi^\vee(k\Lambda_0) = \Pi^\vee_q = \{\al^\vee_1, \ldots, \al^\vee_\ell, -\theta^\vee + q K \}.
\]
And in turn we must have
\[
k = -h^\vee + p/q, \qquad p \geq h^\vee.
\]
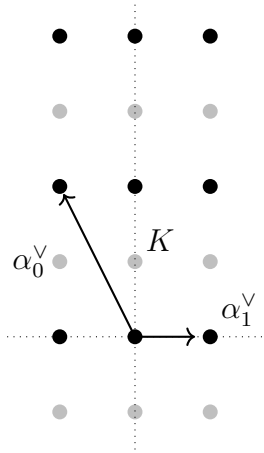
\begin{figure}[H]
    \centering
\begin{tikzpicture}
  % Grid of dots
  \foreach \x in {-1,0,1}
    \foreach \y in {0,2,4}
      \node[circle,fill,inner sep=2pt] at (\x,\y) {};

  \foreach \x in {-1,0,1}
    \foreach \y in {-1,1,3}
      \node[circle,fill=lightgray,inner sep=2pt] at (\x,\y) {};

  % Dotted axes
  \draw[dotted] (-1.7,0) -- (1.7,0);
  \draw[dotted] (0,-1.5) -- (0,4.5);

  % Thick arrows
  \draw[->,thick] (0,0) -- (0.8,0);
  \draw[->,thick] (0,0) -- (-0.95,1.9); 

  % Label (0,1)
  \node[left] at (-1, 1) {$\alpha_0^\vee$};
  \node[above right] at (1,0) {$\alpha_1^\vee$};
  \node[above right] at (0,1) {$K$};
\end{tikzpicture}

    \caption{The coroot subsystem $\what\D^\vee(k\Lambda_0)$ for $\what{\mathfrak{sl}}_2$ and $k$ admissible with denominator $q=2$}
    \label{fig:q=3.simp}
\end{figure}

\begin{rem}
If $\g$ is simply laced, then the intersection of $\what{\D}^{\vee, \text{re}}$ with the set $s K + \h$ is, depending on the value of $s$, either empty or of the form $\D^\vee$, where $\D^\vee$ is the coroot system of $\g$. This is no longer the case when $\g$ is not simply laced; in this case the intersection is either $\D^{\vee}$ or its proper subset $\D^{\vee}_{\text{short}}$ consisting of short coroots. There is therefore a marked difference in the structure of $\Pi^\vee_q$ between the cases $(q, r^\vee) = 1$ and $(q, r^\vee) \neq 1$ when $\g$ is non simply laced with lacing number $r^\vee$. This issue is present throughout the whole story, but we ignore it for the sake of simplicity.
\end{rem}

We now state an important generalisation, due to Arakawa, of Theorem \ref{thm:aff.k.Zplus}
\begin{thm}[Arakawa, {\cite{A.rat.O}}]\label{thm:A.rat.O}
Let $k$ be an admissible level with denominator $q$. Recall $V^k(\g)$ the universal affine vertex algebra of level $k$ and $L_k(\g)$ its simple quotient. The $V^k(\g)$-module $L_k(\lambda)$ descends to a $L_k(\g)$-module if and only if $\lambda$ is admissible with $\Pi^\vee(\lambda)$ conjugate to $\Pi^\vee_q$.
\end{thm}

The proof of this theorem uses a number of advanced techniques. The ``only if'' part is proved by reduction to the case of $\g = \mathfrak{sl}_2$ via \emph{semi-infinite reduction}. The case of $\mathfrak{sl}_2$, in turn, had been settled earlier by Adamovic and Milas using the Zhu algebra \cite{AM95}.

The ``if'' part is proved first for $G$-integrable admissible weights, and the case of general admissible weight is reduced to the $G$-integrable case using the theory of characteristic varieties. For $\g$ a finite dimensional semisimple Lie algebra, and a two-sided ideal $I \subset U(\g)$, the \emph{characteristic variety} $V(I) \subset \h^*$ is introduced as the spectrum of the image of $I$ under the projection to the second factor in
\[
U(\g) \rightarrow (\n_- U(\g) + U(\g) \n_+) \oplus U(\h).
\]
If $\lambda \in V(I)$ then $I L(\lambda) = 0$. For reflections $r_i$ in simple roots, it is a classical fact due to Joseph \cite{Joseph.76} that if $\lambda \in V(I)$ then $r_i \circ \lambda \in V(I)$ too. A generalisation to affine case allows the reduction from general admissible weights to $G$-integrable admissible weights to be carried out.

The case of $G$-integrable $L_k(\lambda)$ is handled using a theory of ``Kazhdan-Lusztig translation functors'' developed by Feigin and Malikov. Although it takes us a little far afield to describe this theory, and we make no pretense at completeness, we sketch the ideas in Section \ref{subsec:KL} below.

\section{Characters formulas for irreducible $\what\g$-modules}

We refer to Kac's book {\cite{Kac.IDLA}} for details. Let $k \in \Z_+$ and recall the set
\[
P_+^k = \{\lambda \mid \left<\lambda, \theta^\vee\right> \leq k\}
\]
of dominant integral weights of level $k$. The corresponding irreducible highest weight $\what\g$-modules $L_k(\la)$ are examples of \emph{integrable} $\what\g$-modules, and their formal characters are given by the Weyl-Kac character formula:
\[
\chi_{L(\lambda)} = \sum_{w \in \what{W}} \varepsilon(w) \chi_{M(w \circ \lambda)},
\]
where $\what{W} = W \ltimes t_{Q^\vee}$ is the affine Weyl group, and
\[
\chi_{M(\lambda)} = \frac{e^{\lambda}}{\prod_{\alpha \in \D_+} (1-e^{-\alpha})}
\]
is the character of the Verma module $M(\la)$. Here the notation $w \circ \la$ represents the ``dot action'' of $W$ on the weight $\la \in \what\h^*$, defined to be
\[
w \circ \lambda = w(\lambda + \what\rho) - \what\rho.
\]

As we have indicated in the previous paragraph, the structure of the Weyl group of $\wtil{\g}$ is that of a semidirect product of the finite Weyl group $W$ with an abelian group of ``translations'' $t_\alpha$ indexed by the coroot lattice $\alpha \in Q^\vee$. We make this more explicit now. The actions of $t_\alpha$ on $\h^*$ and $\h$, respectively, are:
\begin{align}\label{eq:transl.k}
\begin{split}
t_\alpha(\lambda) &= \lambda + \lambda(K)\alpha - \left[ (\alpha, \lambda) + \frac{|\alpha|^2}{2} \lambda(K) \right] \delta \\
t_\alpha(x) &= x - \alpha(x) K \bmod{\C d}.
\end{split}
\end{align}

\begin{figure}[H]
\centering
\begin{tikzpicture}[scale=1]
  % Axes
  \draw[-] (-5, 0) -- (5, 0);
  \draw[-] (0, -5) -- (0, 1);
  
  % Points
  \foreach \x/\y in {-4/-4, -3/-2, -1/0, 0/0, 2/-2, 3/-4} {
    \fill (\x, \y) circle (3pt);
  }

  % x Axis Points
  \foreach \x in {4, 3, 2, 1, 0, -1, -2, -3, -4} {
    \fill (\x, 0) circle (1pt);
  }
  
  % y Axis Points
  \foreach \y in {-1, -2, -3, -4} {
    \fill (0, \y) circle (1pt);
  }

  \node at (-0.5, -2) {$-2\delta$};
  \node at (-0.5, -4) {$-4\delta$};
  
  \node at (4, 0.5) {$4\alpha$};
  \node at (2, 0.5) {$2\alpha$};
  \node at (0.3, 0.5) {$0$};
  \node at (-2, 0.5) {$-2\alpha$};
  \node at (-4, 0.5) {$-4\alpha$};

  \node at (0.3, -0.5) {$\Lambda_0$};
  \node at (3.9, -3.7) {$t_\alpha \circ \Lambda_0$};
  \node at (-3.9, -1.7) {$t_{-\alpha} \circ \Lambda_0$};
  
  \node at (-1, -0.5) {$s \circ \Lambda_0$};
  \node at (2.9, -1.7) {$s t_{-\alpha} \circ \Lambda_0$};
  \node at (-3, -4.1) {$st_\alpha \circ \Lambda_0$};
  
  % Parabola
  \draw[dashed, thick] (-4.3, -4.73) parabola bend (-0.5, 0.0833) (3.3, -4.73);
\end{tikzpicture}
\caption{Dot orbit of $\La_0$ under the affine Weyl group}
\label{fig:dot.orbit}
\end{figure}
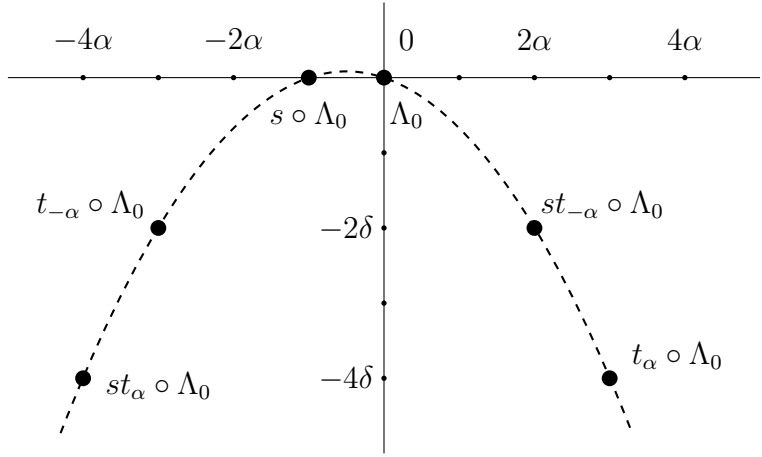

To illustrate we let $\g = \mathfrak{sl}_2$ and we consider the integrable $\what\g$-module $L_1(\g) = L(\Lambda_0)$ of level $k = 1$. We have $h^\vee = 2$ and so $\what\rho = 2\Lambda_0 + \alpha/2$, where $\alpha$ is the unique simple root of $\mathfrak{sl}_2$. A few elements of $\what{W} \circ \Lambda_0$ are shown in Figure \ref{fig:dot.orbit}.
Putting
\[
q = e^{-\delta}, \quad \text{and} \quad y = e^{\alpha/2},
\]
we easily obtain from the Weyl-Kac character formula
\begin{align}\label{eq:L1.sl2.char}
\chi_{L(\Lambda_0)} = \frac{\sum_{n \in \Z} y^{3n} q^{3n^2+n} - \sum_{n \in \Z} y^{3n-1} q^{3n^2-n} }{\prod_{n=1}^\infty (1 - y^{-1} q^{n-1}) (1 - q^n) (1 - y q^n)}.
\end{align}
Notice that the vertex algebra $L_1(\mathfrak{sl}_2)$, as an $\what{\mathfrak{sl}}_2$-module, is isomorphic to $L(\Lambda_0)$. It is not difficult to confirm that under this identification $d = -L_0$. So the graded dimension of the vertex algebra relative to its conformal $\Z_+$-grading, is just the $y=1$ specialisation of \eqref{eq:L1.sl2.char}.

We now make a brief tangential remark. If $Q$ is an even integral lattice, then there is the well known construction of a lattice vertex algebra $V_Q$ \cite{Kac.VA}. This vertex algebra is conformal, and from its construction the graded dimension $\chi_{V_Q}$ is easy to write down, it is
\[
\prod_{n=1}^\infty (1-q^n)^{-1} \cdot \sum_{m \in \Z} e^{\alpha} q^{(\alpha, \alpha)/2}.
\]
The Frenkel-Kac-Segal theorem asserts that if $\g$ is a simply laced simple Lie algebra with root lattice $Q$, then
\[
L_1(\g) \cong V_Q.
\]
So as a consequence of the isomorphism in the case of $\g = \mathfrak{sl}_2$, we ought to recover the identity
\begin{align*}
\prod_{n=1}^\infty (1 - y^{-1} q^{n-1}) (1 - y q^n) \cdot \sum_{m \in \Z} y^m q^{m^2} = \sum_{n \in \Z} y^{3n} q^{3n^2+n} - \sum_{n \in \Z} y^{3n-1} q^{3n^2-n}.
\end{align*}
This is in fact a classical identity.

One of the key features of the admissible weights $\la$ is that the corresponding highest weight modules $L(\la)$, although not integrable, have characters given by formulas essentially similar to the Weyl-Kac character formula. Indeed let $\la$ be an admissible weight for $\what\g$. In \cite{KW88} Kac and Wakimoto proved that
\begin{align}\label{eq:KW.char}
\chi_{L(\lambda)} = \sum_{w \in \what{W}(\lambda)} \varepsilon(w) \chi_{M(w \circ \lambda)},
\end{align}
where now $\what{W}(\lambda) \subset \what{W}$ is the subgroup generated by reflections in the coroots that make up $\Pi^\vee(\lambda)$.

In order to sketch the proof, we require the notions of blocks and linkage. A good introductory reference, for the case of finite dimensional $\g$, is \cite{Humphreys}. See the article \cite{KK} and books \cite{Kac.IDLA} and \cite{Moody.Pianzola} for the Kac-Moody case.

The block decomposition is a fundamental tool for organising our understanding of category $\OO$ of a semisimple (or Kac-Moody) Lie algebra. Basic questions in the theory are to determine (1) for which weights $\mu$ the irreducible module $L(\mu)$ may occur in a composition series of the Verma module $M(\la)$, and (2) with what multiplicities. A definitive answer to the second question came with the Kazhdan-Lusztig conjecture and its solution, and extensions thereof. As for the first question, an important necessary condition (besides the obvious $\la - \mu \in Q_+$) is that there exist $w \in W$ such that
\[
\mu = w(\la+\rho)-\rho.
\]
We write $w \circ \la$ for the term on the right hand side, and say that weights satisfying such a relation are ``linked''. In fact if $\la \notin P$ then not every $w \circ \la$ differs from $\la$ by an element of $Q$, so it is possible to make an obvious refinement to the notion of linkage: we define $W(\lambda) \subset W$ the subgroup generated by reflections in those coroots $\al$ such that $\left<\la +\rho, \al^\vee\right> \in \Z$ (for an admissible weight $\la$ of $\what\g$, this group is precisely the group $\what{W}(\la)$ appearing in the character formula \eqref{eq:KW.char}), and linkage is defined as the equivalence relation generated by $\mu \in W(\la) \circ \la$.

The objects of $\OO$ possessing a composition series all of whose irreducible factors have highest weight in a fixed ``linkage class'', forms an abelian subcategory of $\OO$. There are no nonzero morphisms between objects of these subcategories, which are known as ``blocks'', and $\OO$ decomposes into a direct sum of blocks.

It can be convenient to translate results from one block to another by means of so called translation functors. Let $\OO_\la$ and $\OO_\mu$ denote blocks associated with weights $\la$ and $\mu$. If $\la-\mu \in P$ then the translation functor $T_\la^\mu : \OO_\la \mapsto \OO_\mu$ is defined by
\[
T_\lambda^\mu(M) = \text{pr}_\mu \left( M \otimes L(\nu) \right),
\]
where $\nu$ is the unique representative of $P_+ \cap W(\la-\mu)$. Here $\text{pr}_\mu$ means projection from $\OO$ to the summand $\OO_\mu$. The module $L(\nu)$ is integrable, which turns out to be important to guarantee exactness of the translation functor.

If $\mu$ is in the same Weyl chamber as $\lambda$, or a facet of it (in the sense that for each $\alpha^\vee \in \D^\vee(\lambda)$ we have $\left<\lambda+\rho, \al^\vee\right> > 0$ (resp. $< 0$) implies $\left<\mu+\rho, \al^\vee\right> \geq 0$ (resp. $\leq 0$)), then 
\[
T_\lambda^\mu M(w \circ \lambda) = M(w \circ \mu).
\]
The proof of this fact relies on the following standard lemma (see \cite[Lemma 7.5]{Humphreys}).
\begin{lemma}\label{lem:Hum.comb}
Let $\lambda, \mu \in P_+$, and let ${\nu}$ be the unique element of the intersection of $W(\mu-\lambda)$ with $P_+$. Then the only weight $\lambda+\nu'$ (where $\nu'$ is a weight of $L({\nu})$) that is linked to $\mu$, is $\mu = \lambda+\nu$ itself.
\end{lemma}

Returning now to the character of admissible $\what\g$-modules. By the linkage principle
\begin{align*}
\chi_{L(\lambda)} = \sum_{w \in \what{W}(\lambda)} m(\lambda, w) \, \chi_{M(w \circ \lambda)}
\end{align*}
for some set of integer coefficients $m(\lambda, w)$. Let $s = r_{\alpha_i}$ be a simple reflection in $\what{W}(\lambda)$, and choose $\mu$ so that $\what{W}(\mu-\lambda)$ meets $P_+$, and so that $\left<\mu+\rho, \al^\vee\right> \in \Z_{\geq 0}$ for all $\alpha^\vee \in \Pi^\vee(\lambda)$, with $\left<\mu+\rho, \al^\vee\right> = 0$ for exactly the one coroot $\alpha^\vee = \alpha_i^\vee$. We apply a translation functor to the surjection
\[
M(\lambda) / M(s \circ \lambda) \rightarrow L(\lambda).
\] 
Since $s \circ \mu = \mu$, we get $T_\lambda^\mu L(\lambda) = 0$ (an exact functor sends surjections to surjections), hence
\[
\sum_{w \in \what{W}(\lambda)} m(\lambda, w) \chi_{M(w \circ \mu)} = 0,
\]
and so $m(\lambda, w) + m(\lambda, s \circ w) = 0$. This goes for each $s = r_{\alpha_i}$, and so $m(\lambda, w) = \varepsilon(w)$ for all $w \in \what{W}(\lambda)$, thus proving the formula.

\section{The Kazhdan-Lusztig tensor product}\label{subsec:KL}

The tensor product, as vector spaces, of a $\what\g$-module of level $k$ with one of level $\ell$ is again a $\what\g$-module in the usual way, and has level $k+\ell$. On the other hand Kazhdan and Lusztig have defined a tensor product for modules of a single level $k$, which we describe briefly in this section. The main reference is the original series of four articles \cite{KL1, KL2, KL3, KL4}.

Above we have discussed $G$-integrable modules already, it is convenient to introduce notation for this now. Let $\KL_k$ denote the category of $G$-integrable positive energy $\what\g$-modules of level $k$. This category contains the Weyl modules $V_k(\lambda)$ for $\lambda \in P_+$ in particular.

Let $\Gamma(\C P^1 \backslash \{0, 1, \infty\}, \OO)$ denote the ring of rational functions on the complex projective space $\C P^1$ of dimension one, with poles only at $0$, $1$ and $\infty$. Let $t$ denote the standard coordinate on $\C = \C P^1 \backslash \{\infty\}$. Following Kazhdan and Lusztig, we consider the Lie algebras
\begin{align*}
\Gamma &= \g \otimes \Gamma(\C P^1 \backslash \{0, 1, \infty\}, \OO) \\
\text{and} \quad \what\Gamma &= \Gamma \oplus \C K
\end{align*}
with Lie brackets
\begin{align*}
[x \, f, y \, g] &= [x, y] \, fg \\
\text{and} \quad [x \, f, y \, g] &= [x, y] \, fg + \res_{t=0} f'(t) g(t) K.
\end{align*}
There is now a homomorphism $\what\Gamma \rightarrow \what\g$ induced by Laurent series expansion of rational functions
\[
i_0 : \Gamma(\C P^1 \backslash \{0, 1, \infty\}, \OO) \rightarrow \C((t))
\]
at $0$ in powers of $t$.

Let $A, B \in \KL_k$. We turn the vector space tensor product
\[
W = A \otimes B
\]
into a $\what\Gamma$-module, by equipping it with the action defined by
\[
x f (a \otimes b) = (x \, i_1 f) a \otimes b + a \otimes (x \, i_\infty f) b, \qquad K = -k \, \text{Id}.
\]
The choice of action of $K= -k \, \text{Id}$ here was not made arbitrarily, it is the only choice consistent with the action of terms $x f$ as defined. This is a consequence of the residue theorem. We denote by $Z$ the algebraic dual of $W$. This carries an action of $\what\Gamma$ of level $k$.

Now $Z$ is almost the module we want; a module of level $k$ constructed from a pair of modules $A$ and $B$ of level $k$. However, as an algebraic dual of an infinite dimensional vector space, $Z$ has larger cardinality than $A$ and $B$, whereas we ultimately want to construct a tensor product $\dot\otimes : \KL_k \times \KL_k \rightarrow \KL_k$.

Inside $Z$ we consider subspaces
\[
Z^N = \{\varphi \in Z \mid \varphi(\Gamma_0^N W) = 0 \},
\]
where $\Gamma_0 \subset \what\Gamma$ denotes the subspace of $\g$-valued functions that vanish at $0$. These form an increasing sequence, and we denote by $Z^\infty$ their union. Let $f(t) \in \C((t))$ and $x \in \g$. How do we define the action of $x \, f$ on $\varphi \in Z^\infty$? Noting that by definition $\varphi \in Z^N$ for some $N \in \Z_+$, the answer is to approximate $f(t)$ by a rational function $\wtil{f}(t) \in \G(\C P^1 \backslash \{0, 1, \infty\}, \OO)$, and define $x\, f \cdot \varphi$ by
\[
[x\, f \cdot \varphi](a \otimes b) = \varphi(x\, \wtil{f} (a \otimes b)).
\]
If the chosen approximation is accurate to order at least $t^N$, then $\varphi$ ought to annihilate the difference between $x\, f$ and $x\, \wtil f$, and so this definition of the action of $\what\g$ is reasonable.
\begin{thm}[KL1, Section 7]
If $A, B \in \KL_k$ then $Z^\infty$ is again a $\what\g$-module in $\KL_k$.
\end{thm}

If we are given three objects $A$, $B$ and $C$ of $\KL_k$ then we can consider the action of $\G$ on $A \otimes B \otimes C$ given by
\[
x \, f (a \otimes b \otimes c) = (x \, i_1 f) a \otimes b \otimes c + a \otimes (x \, i_\infty f) b \otimes c + a \otimes b \otimes (x \, i_0 f) c.
\]
This is well-defined because of the residue theorem. Then we consider the vector space
\begin{align}\label{eq:conf.block.def}
\CC(A, B, C) = \left(\frac{A \otimes B \otimes C}{\Gamma \cdot (A \otimes B \otimes C)}\right)^*,
\end{align}
which is finite dimensional under favourable circumstances, and is known as the space of conformal blocks. From $A$ and $B$ we construct $Z^\infty$ as above, and then it turns out that
\[
\Hom_{\what\g}(C, Z^\infty) \cong \CC(A, B, C).
\]
Reflecting on the tensor-hom adjunction, this relation shows that $Z^\infty$ behaves like the dual of a tensor product. This is not surprising since we have defined $Z^\infty$ inside an algebraic dual. The final step in defining $A \dot\otimes B$ is to take a graded (or ``contragredient'') dual $D(-)$, twisting the action of $\what\g$ in the process. We do not enter into details here, merely remarking that the procedure is much the same as that employed to define duals in category $\OO$ of a finite dimensional semisimple Lie algebra \cite{Humphreys}. The Kazhdan-Lusztig tensor product is thus the contragredient dual of $Z^\infty$:
\[
A \dot\otimes B = D(Z^\infty).
\]

The construction of $A \dot\otimes B$ can be carried out with $A$ and $B$ located, not at $\infty$ and $1$ respectively, but at any pair of points $(p_1, p_2)$ distinct from $0$ and from each other. An important part of the theory, which it would be too large a detour to describe in detail in these notes, is the construction of a connection on the vector bundle $\CC$ over the configuration space of such pairs $(p_1, p_2)$, whose fibres are the corresponding spaces of conformal blocks. Using this connection, one can let the points move and, importantly, loop around one another and exchange positions. In this way for example, braiding isomorphisms
\[
A \dot\otimes B \rightarrow B \dot\otimes A
\]
can be defined. Associativity isomorphisms are also defined, and correspond to configurations of three points (distinct from $0$) passing to degenerate configurations, in which two points ``collide'' in distinct ways. In the end a braided monoidal category is obtained. For a complete definition, and theory, of (braided) monoidal category, see \cite{EGNO}. For an account emphasising links with vertex algebras, see \cite{BK}.
\begin{defn}
A monoidal category is a category $\CC$ with a bifunctor $\otimes : \CC \times \CC \rightarrow \CC$ and natural isomorphism
\[
a_{X, Y, Z} : (X \otimes Y) \otimes Z \rightarrow X \otimes (Y \otimes Z),
\]
a unit $\mathbf{1} \in \CC$ and unit natural isomorphism, satisfying certain identities known as the pentagon and triangle identities. A braiding on a monoidal category is a further natural isomorphism
\[
c_{X, Y} : X \otimes Y \rightarrow Y \otimes X,
\]
compatible with the associativity and unit isomorphisms through identities known as the hexagon identities.
\end{defn}

The Kazhdan-Lusztig category is easier to describe for irrational level $k$ than for rational level. In this case the Weyl modules $V_k(\lambda)$, for $\la \in P_+$, are irreducible, and the Kazhdan-Lusztig tensor product is compatible with induction in the sense that
\[
V_k(\lambda) \dot\otimes V_k(\mu) \cong \bigoplus_{\nu \in P_+} V_k(\nu)^{\oplus m_\nu}, \quad \text{where} \quad L(\lambda) \otimes L(\mu) \cong \bigoplus_{\nu \in P_+} L(\nu)^{\oplus m_\nu}.
\]
The situation at rational values of the level $k$ is more complicated, but what can be said in general is that
\[
V_k(\lambda) \dot\otimes V_k(\mu)
\]
possesses a filtration, whose composition factors are the Weyl modules $V_k(\nu)$ for $\nu \in P_+$ with the same multiplicities $m_\nu$ as for irrational level $k$ \cite{KL4} \cite[Theorem 3.1.2]{FM}.

The category $\KL_k$ is a subcategory of $\OO$ and as such splits into blocks, corresponding to the equivalence relation
\[
\mu = w \circ \lambda, \qquad w \in \what{W}_k.
\]
Here the group $\what{W}_k = W \ltimes p Q^\vee = \what{W}(k\Lambda_0)$ appears naturally since objects of $\KL_k$ are by definition $G$-integrable, and by \eqref{eq:transl.k} translation $t_\al(\la)$ may satisfy $t_\al(\la)-\la \in Q$ only if $t_\al(\la)-\la \in pQ$.

The geometry of weights under the action of $W_k$ is substantially similar to the case of $\g$ semisimple, and linkage and translation functors can be defined in this context. For $\lambda, \mu \in P_+^k$ regular, let $\nu$ be the unique weight in $P_+^k \cap W(\mu-\lambda)$, then for $M \in \OO_\la$ one defines
\[
T_\lambda^\mu(M) = \text{pr}_\mu(M \dot\otimes V_k(\nu)).
\]
Notice that the tensor product is taken with the Weyl module $V_k(\nu)$ rather than the irreducible quotient $L(\la)$. It can be proved that $T_\la^\mu$ sends Weyl modules to Weyl modules, and more specifically that
\[
T_\lambda^\mu(V_k(w \circ \lambda)) \cong V_k(w \circ \mu),
\]
whenever $w \circ \lambda, w \circ \mu \in P_+$. This is a consequence of a version of Lemma \ref{lem:Hum.comb} for $\what{W}_k$, which is contained in Jantzen's book \cite[Lemma 7.7]{Jantz.alg}, see also \cite[Lemma 2.1.1]{FM}. Next we would like to prove that
\begin{align}\label{eq:KL.trans}
T_\lambda^\mu(L_k(w \circ \lambda)) \cong L_k(w \circ \mu).
\end{align}
Conceptually this is done using the notions of duality and rigidity in monoidal categories.

\begin{defn}
For objects $X$ and $Y$ in a monoidal category, a duality datum (exhibiting $X$ as a left dual of $Y$ and $Y$ as a right dual of $X$) is a pair of morphisms
\[
\ev : Y \otimes X \rightarrow \mathbf{1}, \qquad \coev : \mathbf{1} \rightarrow X \otimes Y
\]
such that the compositions
\[
X \rightarrow \mathbf{1} \otimes X \rightarrow X \otimes Y \otimes X \rightarrow X \otimes \mathbf{1} \rightarrow X 
\]
and
\[
Y \rightarrow Y \otimes \mathbf{1} \rightarrow Y \otimes X \otimes Y \rightarrow \mathbf{1} \otimes Y \rightarrow Y 
\]
are the identities.
\end{defn}
These axioms capture the essential properties of the notion of dual of a finite dimensional vector space over a field $\F$ (and of dual module for modules of a group or Lie algebra, etc.). The evaluation and coevaluation morphisms for such a vector space $V$ and its dual $V^*$ are
\begin{align*}
\ev &: V^* \otimes V \rightarrow \F, \qquad \varphi \otimes v \mapsto \varphi(v) \\
\text{and} \quad
\coev &: \F \rightarrow V \otimes V^*, \qquad 1 \mapsto \sum_i v_i \otimes \varphi_i,
\end{align*}
where $\{v_i\}$ is any basis of $V$ and $\{\varphi_i\}$ is its dual basis of $V^*$.

The following is a standard result in the theory.
\begin{prop}
If $(X, X^*)$ is a duality datum in an abelian monoidal category, i.e., if there are evaluation and coevaluation maps making $X^*$ a left dual of $X$, then
\[
\Hom(U, X \otimes V) \cong \Hom(X^* \otimes U, V)
\]
and $X \otimes (-)$ is an exact functor.
\end{prop}
A monoidal category $\CC$ is said to be \emph{rigid} if there is a contravariant functor $X \mapsto X^*$ and natural duality data relating $X$ and $X^*$ for all $X \in \CC$.

The category $(\KL_k, \dot\otimes)$ is rigid if $k \notin \Q$, while for $k \in \Q_{>0}$ there are in general some objects that do not admit a duality datum (highest weight modules with highest weight that lies on an affine wall for $\what{W}_k$). In this case Frenkel-Malikov \cite{FM} observe that $V_k(\la)$ possesses a duality datum for regular weights $\lambda$, with dual $V_k(\la)^* = V_k(-w_\circ(\la))$, where $w_\circ$ is the longest element in the Weyl group of $\g$.

In general if $V_k(\lambda)$ fits into a duality datum, we have exactness of $T_\la^\mu$, then a standard argument (see {\cite[Proposition 7.8]{Humphreys}} for the simplest iteration of this argument, and \cite{FM} for the present context) implies \eqref{eq:KL.trans}.

Let $\la \in P_+$ and recall that the irreducible $\what\g$-module $L_k(\la)$ is a quotient of the Weyl module $V_k(\la)$. In particular $L_k(\g) = V_k(0) / I_k$ where $I_k$ is both a $\what\g$-submodule and a vertex algebra ideal. Suppose now that $\la \in P_+^k$ is a regular dominant integral weight. In order to prove that $L_k(\la)$ is a $L_k(\g)$-module, we should like to prove that
\begin{align}\label{eq:translation.desire}
L_k(\lambda) \cong V_k(\lambda) / I_k \cdot V_k(\lambda).
\end{align}

Let $V$ be a vertex algebra, $V^k(\g)$ in the case at hand, let $M_1$, $M_2$ and $N$ be three $V$-modules, and let $I \subset V$ be an ideal. The following is an easy but key observation:
\begin{align*}
\CC(M_1/IM_1, M_2, N) = \CC(M_1, M_2/IM_2, N).
\end{align*}
Indeed let us write $\Gamma$ for global sections as in \eqref{eq:conf.block.def}. Then the left and right hand sides are identified, respectively, with the subsets
\begin{align*}
&\{\varphi : M_1 \otimes M_2 \otimes N \rightarrow \C \mid \text{$\varphi(i_{(n)}M_1, M_2, N) = 0$ for all $i \in I$ and $n \in \Z$}\} \\
\text{and} \quad &\{\varphi : M_1 \otimes M_2 \otimes N \rightarrow \C \mid \text{$\varphi(M_1, i_{(n)}M_2, N) = 0$ for all $i \in I$ and $n \in \Z$}\} 
\end{align*}
of
\begin{align*}
\CC(M_1, M_2, N) = \{\varphi : M_1 \otimes M_2 \otimes N \rightarrow \C \mid \text{$\varphi(\Gamma \cdot (M_1 \otimes M_2 \otimes N)) = 0$}\}.
\end{align*}
All we have to do is convince ourselves that these two subsets coincide.

From \eqref{eq:KL.trans} we deduce $L_k(\la) \cong T_0^\lambda(L_k(\g))$. We combine this with the other ingredients above in the following simple calculation:
\begin{align*}
\Hom_{\what\g}(L_k(\la), D(N)) &\cong 
\Hom_{\what\g}(T_0^\lambda(V_k(0) / I_k), D(N)) \\
&\cong \CC(V_k(0) / I_k, V_k(\lambda), N) \\
&= \CC(V_k(0), V_k(\lambda) / I_k \cdot V_k(\lambda), N) \\
&\cong \Hom_{\what\g}(V_k(\lambda) / I_k \cdot V_k(\lambda), D(N)).
\end{align*}
So by the representing property of the tensor product, we have \eqref{eq:translation.desire} as desired.

\section{The BRST complex and its quantisation}

Hamiltonian reduction of symplectic varieties admitting Lie group symmetry is an important basic construction in mathematical physics, admitting many refinements and variations. Firstly the notion of Marsden-Weinstein quotient of a symplectic manifold by a good Hamiltonian group action is given a homological interpretation in the form of the BRST complex \cite{KS}. This homological construction opens the door both to quantisation, and to generalisation to the infinite dimensional context, where it yields a Hamiltonian reduction functor for vertex algebras. Via this construction the \emph{affine $W$-algebras} are defined and studied.

Let $M$ be a smooth Poisson manifold equipped with a Hamiltonian action of a Lie group $N$. Let $\n$ be the Lie algebra of $N$. Then $\OO = C^\infty(M)$ is a Poisson algebra, and the Hamiltonian action basically amounts to a homomorphism of Lie algebras
\[
\alpha : \n \rightarrow \OO.
\]
The moment map $\mu : M \rightarrow \n^*$ takes a point $p \in M$ to the linear functional on $\n$ that sends $x \in \n$ to the number obtained by evaluating $\alpha(x) \in \OO$ at $p$.

The algebra of functions on $\mu^{-1}(0)$ is the quotient $\OO / (\alpha(\n))$, where $(\ldots)$ denotes ``ideal generated'' in the commutative algebra $\OO$. Under favourable hypotheses on the structure of the action, this quotient algebra is resolved by the Koszul complex:
\[
 \OO \otimes \wedge^\bullet\n,
\]
the differential being given by
\begin{align*}
f \otimes \phi &\mapsto \alpha(\phi) f, \\
f \otimes \phi_1 \wedge \phi_2 &\mapsto \alpha(\phi_1) f \otimes \phi_2 - \alpha(\phi_2) f \otimes \phi_1,
\end{align*}
etc.

So far we have not used either the Lie algebra structure of $\n$ nor the Poisson bracket in $\OO$. We do use these structures, however, to furnish the whole Koszul complex with the structure of a complex of $\n$-modules. Finally the subalgebra of $\OO/(\alpha(\n))$ consisting of $N$-invariant functions on $\mu^{-1}(0)$ is resolved by the Chevalley-Eilenberg complex of the $\n$-module $\OO \otimes $, namely $(\OO \otimes \wedge^\bullet\n) \otimes \wedge^\bullet\n^*$. The map $\OO \rightarrow \OO \otimes \n^*$ is the natural one obtained from the action of $\n$ through tensor-hom adjunction:
\[
\Hom(\n \otimes \OO, \OO) \cong \Hom(\OO, \OO \otimes \n^*).
\]
More explicitly, if $\{\phi_1, \ldots, \phi_n\}$ is a basis of $\n$ and $\{\phi_1^*, \ldots, \phi_n^*\}$ is the corresponding dual basis of $\n^*$, then
\[
f \mapsto \sum_{i=1}^n \{\alpha(\phi_i), f\} \otimes \phi_i^*.
\]
The other vertical morphisms are constructed in a similar way.
\begin{align*}
\xymatrix{
& & & \\
{} \ar@{->}[r] & (\OO \otimes \wedge^2\n) \otimes \n^* \ar@{->}[r] \ar@{->}[u] & (\OO \otimes \n) \otimes \n^* \ar@{->}[r] \ar@{->}[u] & \OO \otimes \n^* \ar@{->}[u] \\
{} \ar@{->}[r] & \OO \otimes \wedge^2\n \ar@{->}[r] \ar@{->}[u] & \OO \otimes \n \ar@{->}[r] \ar@{->}[u] & \OO \ar@{->}[u] \\
}
\end{align*}

So we have constructed a bicomplex, from which we may pass to the total complex $C(\n, \OO)$ with total differential $d_{\text{Tot}}$. Under favourable circumstances, the total complex has cohomology concentrated in degree $0$, recovering the ring of $N$-invariant functions on $\mu^{-1}(0)$.

Now the quantisation proposal: Suppose we knew how to quantise the Poisson manifold $M$, that is, suppose we had a filtered associative algebra $A$ whose associated graded (i.e., classical limit) recovers $\OO = \gr(A)$. Suppose further that we knew how to quantise $M$ equivariantly, that is, suppose $A$ comes with a compatible homomorphism $\alpha : U(\n) \rightarrow A$. Then a natural quantisation can be found for the BRST complex, and thus for $\mu^{-1}(0)/N$.

Indeed, let $\ma$ be a vector space and $\left<\cdot,\cdot \right> : \ma \times \ma \rightarrow \C$ a symmetric bilinear form. The corresponding Clifford algebra is the quotient
\[
\Cl(\ma) = \frac{T(\ma)}{(ab + ba - \left<a, b\right>)}
\]
of the tensor algebra $T(\ma)$. We may also consider $\ma$ as a purely odd vector superspace, and $T(\ma)$ as an associative superalgebra. The canonical grading on $T(\ma)$ by word length induces a filtration on $\Cl(\ma)$, and the associated graded is isomorphic to the exterior algebra $\wedge(\ma)$. Let $A$ be a filtered associative algebra, we can form the filtered associative superalgebra
\[
A \otimes \Cl(\ma).
\]

Now we suppose, in the context of the Hamiltonian action of $N$ on $M$ as above, that our associative algebra $A$ is a quantisation of $M$ in the sense that the associated graded $\gr(A)$ recovers the Poisson algebra $\OO$. If we take $\ma = \n \oplus \n^*$, where $\n$ is our Lie algebra from above, then $\ma$ can be given the tautological symmetric bilinear form $\left<\cdot, \cdot\right>$ of dual pairing between $\n$ and $\n^*$. Also $\ma$ acquires a semidirect product Lie algebra structure, with Lie brackets as follows given by the adjoint and coadjoint actions
\begin{align*}
\n \times \n &\rightarrow \n \\
\n \times \n^* &\rightarrow \n^* \\
\n^* \times \n^* &\rightarrow 0.
\end{align*}
The form $\left<\cdot, \cdot\right>$ is invariant for this Lie bracket. The associated graded of the filtered superalgebra $A \otimes \Cl(\ma)$ now recovers
\[
\OO \otimes \wedge^\bullet(\n \oplus \n^*).
\]

Somewhat remarkably, there exists in $A \otimes \Cl(\ma)$ a cubic element
\[
Q = \sum_i \alpha(\phi_i) \otimes \phi_i^* - 1 \otimes \frac{1}{2} \sum_{i,j} \phi_i^* \phi_j^* [\phi_i, \phi_j]
\]
(the summation indices run over a basis of $\n$) with the property that $Q^2 = 0$ (which implies that $\ad(Q)^2 = 0$ and hence $\ad(Q)$ defines a differential in $A \otimes \Cl(\ma)$), and in the associated graded we recover the total differential as
\[
d_{\text{Tot}} = \{Q, -\}.
\]
The proposal, therefore, is to regard
\[
H^\bullet(A \otimes \Cl(\ma), \ad(Q))
\]
as a quantisation of $\mu^{-1}(0) / N$.

For simplicity we have described the reduction of $\mu^{-1}(0)$ and its quantisation. All of this can be generalised to $\mu^{-1}(p)$ for suitable $p \in M$ by modifying the differential $Q$. 

\section{Slodowy slices and Drinfeld-Sokolov reduction}

We wish to apply the machinery described above to a particular class of examples of Hamiltonian reductions, known as Slodowy slices. In this particular case, the relevant structures all have vertex algebraic analogues, and so a parallel story can be developed for vertex algebras. In this context the construction is known as quantised Drinfeld-Sokolov reduction.

Let $\g$ be a simple Lie algebra with invariant bilinear form $(\cdot, \cdot)$, and let $f \in \g$ be a nilpotent element. The Jacobson-Morozov theorem \cite{Coll.McGov} asserts that $f$ can be included into an $\mathfrak{sl}_2$-triple $\{e, h, f\} \subset \g$. The endomorphism $\ad(h)$ induces a decomposition of $\g$ into eigenspaces
\begin{align}\label{eq:Dynkin.grading}
\g = \bigoplus_{j \in \Z} \g_j, \qquad \ad(h)|_{\g_j} = j.
\end{align}
At the same time, $\g$ decomposes into finite dimensional irreducible $\ma$-modules, where $\ma \subset \g$ denotes the subalgebra isomorphic to $\mathfrak{sl}_2$ spanned by $e$, $h$ and $f$. If $\g_j$ vanishes for odd $j$ then we say the nilpotent element $f$ is even. Although it is not essential, in these notes we restrict our discussion to even nilpotent elements.

\begin{figure}[H]
\centering
\begin{tikzpicture}
  \begin{scope}[opacity=0.4]
    \fill[blue!30] (-0.47,0) -- (0,0.47) -- (1.13,-0.67) -- (5.33,-0.67) -- (5.33,-1.33) -- (0.87,-1.33) -- cycle;
  \end{scope}

  \begin{scope}[opacity=0.4]
    \fill[blue!30] (3,-3.47) -- (1.2,-1.67) -- (4.8,-1.67) -- cycle;
  \end{scope}
  
  % Dots and labels
  \node[circle,fill,inner sep=2pt] (dot1) at (0,0) {};
  \node[circle,fill,inner sep=2pt] (dot2) at (1,-1) {};
  \node[circle,fill,inner sep=2pt,label={$h$}] (dot3) at (2,-2) {};
  \node[circle,fill,inner sep=2pt] (dot4) at (3,-3) {};
  \node[circle,fill,inner sep=2pt,label={$e$}] (dot5) at (3,-1) {};
  \node[circle,fill,inner sep=2pt] (dot6) at (4,-2) {};
  \node[circle,fill,inner sep=2pt] (dot7) at (5,-1) {};
  \node[circle,fill,inner sep=2pt,label={$f$}] (dot8) at (1,-3) {};
  \node[circle,fill,inner sep=2pt] (dot9) at (2,-4) {};
  \node[circle,fill,inner sep=2pt] (dot10) at (1,-5) {};

  % Dotted lines
  \draw[dashed] (0,1) -- (4,-3);
  \draw[dashed] (-1,0) -- (3,-4);

  \node at (0,2) {$\n$};
  \node at (-1,1) {$\g_0$};
  \node at (-2,0) {$\n_-$};

  \node at (5,-2.5) {$[f, \n]$};
  \node at (4,-0.2) {$\g^e$};
\end{tikzpicture}
\caption{Decomposition of $\g$ with $\g^e$ and $[f, \n]$ highlighted.}
\label{fig:nilp.decomp}
\end{figure}
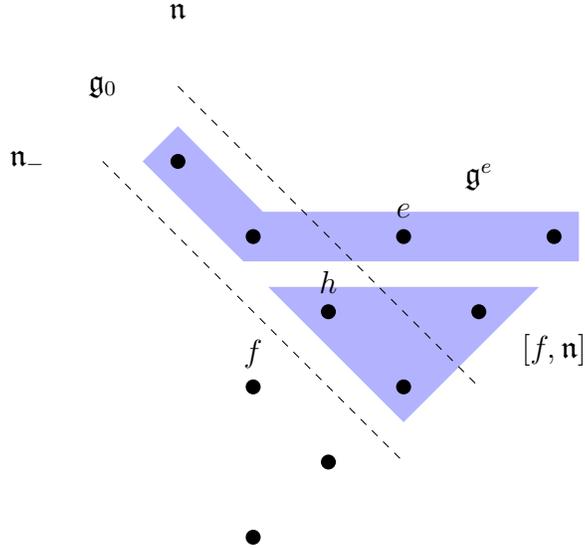
Let us consider the subalgebras $\n = \g_{>0}$ and $\g^e = \{x \in \g \mid [e, x] = 0\}$, and the subspace $\n^\perp = \{x \in \g \mid (x, \n) = 0\}$ (which in this case is also a subalgebra, in fact). It is not difficult to see that
\begin{align}\label{eq:inf.GG}
[\n, f] + \g^e = \g_{\geq 0} = \n^\perp.
\end{align}
This equality has been upgraded by Gan and Ginzburg \cite{GG2002} to an isomorphism of algebraic (Poisson) varieties
\[
N \times (f + \g^e) \cong f + \n^\perp,
\]
where $N$ is the unipotent Lie group with Lie algebra $\n$ and the isomorphism is given by the action map. The isomorphism \eqref{eq:inf.GG} is recovered as the infinitesimal expression. This motivates us to interpret the following natural maps
\[
\alpha : \n \rightarrow \g \subset S(\g) = \C[\g^*], \qquad \mu : \g^* \rightarrow \n^*,
\]
as Hamiltonian action and momentum map. The affine variety
\[
S = S_f = f + \g^e
\]
is known as the \emph{Slodowy slice} associated with $f$. The theorem of Gan and Ginzburg, generalising the case of principal nilpotent due to Kostant, is as follows.
\begin{thm}[\cite{GG2002}]
Let $\chi \in \g^*$ identified with $f \in \g$ under the invariant bilinear form $(\cdot, \cdot)$. Then the Hamiltonian reduction
\[
\mu^{-1}(\chi)/N
\]
is an affine Poisson variety, naturally isomorphic to the Slodowy slice
\[
S = f + \g^e.
\]
\end{thm}

In the case of the Poisson varieties $\g^*$ and $S_f$, enough of the notions and objects have clear infinite dimensional analogues to permit the rendition of the BRST construction to the context of vertex algebras.

We take the universal affine vertex algebra $V = V^k(\g)$ and the charged free fermions vertex superalgebra $F = F^{\text{ch}}(\n)$ based on the subspace $\n \subset \g$. Recall that $F$ is generated by odd fields $\phi_i$ and $\phi_i^*$ indexed by dual bases of $\n$ and $\n^*$, with
\[
[{\phi_i}_\lambda \phi_j^*] = \delta_{ij}.
\]
We overload notation by writing $\alpha(\phi)$ for the field in $V$ corresponding to $\alpha(\phi) \in \g$. Let $C = V \otimes F$ and put
\[
\Qst = \sum_i \alpha(\phi_i) \otimes \phi_i^* - \mathbf{1} \otimes \frac{1}{2} \sum_{i, j} : [\phi_i, \phi_j] \, \phi_i^* \phi_j^*:.
\]
Using the $\lambda$-bracket calculus, we may compute explicitly
\[
[\Qst_\lambda \Qst] = 0,
\]
and so $d_{\text{st}} = \Qst_{(0)}$, which satisfies $d_{\text{st}}^2 = 0$, will be a differential.

The definition of \emph{semi-infinite homology} has been introduced by B. Feigin \cite{Feigin.1984.semiinf}, and the complex $(C, d_{\text{st}})$ coincides with the complex of semi-infinite homology of the $\n[t, t^{-1}]$-module $V$ {\cite[p. 246]{FBZ}}. For this reason its homology is often met in the literature under the following notation
\[
H_i(C, d_{\text{st}}) = H^{\tfrac{\infty}{2}+i}(\n[t, t^{-1}], V).
\]

Once again, this homology corresponds to a quantisation of a Hamiltonian reduction of $\mu^{-1}(0)$. The modification from $\mu^{-1}(0)$ to $\mu^{-1}(\chi)$ entails replacement of $\Qst$ with
\[
\Qst + p, \quad \text{where} \quad p = \sum_i (f, \phi_i) \mathbf{1} \otimes \phi_i^*.
%
%\sum_i \left( e_i + (f, e_i) \right) \otimes \phi_i^* - \frac{1}{2} \otimes \sum_{i, j} :[\phi_i, \phi_j]\phi_i^* \phi_j^*:
\]
It may be verified, again by direct calculation with $\lambda$-brackets, that $d = (\Qst + p)_{(0)}$ satisfies $d^2 = 0$ {\cite[p. 316]{KRW}}. Finally the homology vertex algebra
\[
W = H(C, d)
\]
is the affine $W$-algebra.

\section{Structure, modules and characters of $W$-algebras}\label{sec:str-mod-char}

We have just defined a vertex algebra, really a class of vertex algebras, $W = H(C, d)$. But we have little idea what these algebras ``look like''. Vertex algebras are infinite dimensional, but a grading (for example a conformal grading) with finite dimensional components allows us to write down a character, or graded dimension, of $V$. Such characters are a useful tool in the theory.

We have already seen that $V = V^k(\mathfrak{sl}_2)$ has a conformal structure, and thus a grading (in this case a $\Z_+$-grading) by eigenvalues of $L_0$. Counting the monomials \eqref{eq:Vkg.monom} with conformal weight yields the character
\begin{align*}
\chi_V(q) = \tr_V q^{L_0} 
&= \prod_{n=1}^\infty \frac{1}{(1- q^n)^3}.
\end{align*}
The vertex superalgebra $F = \Fch(\C \phi)$ is also conformal, with conformal vector
\[
L^{\text{ch}} = :(T\phi) \phi^*: = -:\phi^* (T\phi):.
\]
Under this conformal structure we have $\Delta(\phi) = 0$ and $\D(\phi^*) = 1$. For a vector superspace $U = U_{\ov 0} \oplus U_{\ov 1}$ the super dimension $\dim(U_{\ov 0}) - \dim(U_{\ov 1})$ is more natural than the dimension. The supercharacter of $F$ is
\[
\chi_F(q) = \prod_{n=1}^\infty (1-q^{n-1}) (1-q^n) = 0,
\]
in this case the even and odd components of $F$ have matching characters, which cancel out. An additional $\Z$-grading by ``charge'' can be introduced in $F$ so that $F = \bigoplus_{m \in \Z} F^m$, in which $\phi^*$ has charge $+1$ and $\phi$ has charge $-1$. In fact, a $\lambda$-bracket computation reveals that the charge grading coincides with the grading by eigenvalues of $\xi_0$, where $\xi = :\phi^* \phi:$. We now have a two-variable character
\begin{align}\label{eq:F.char.2var}
\chi_F(q, y) = \tr_F y^{\xi_0} q^{L_0} = \prod_{n=1}^\infty (1 - y^{-1} q^{n-1}) (1 - y q^{n}).
\end{align}

\begin{figure}[H]
\centering
\begin{tikzpicture}
%%%%%%%%
% Vacuum picture
  \draw[-] (0, 0) -- (-3, -3); 
  \draw[-] (0, 0) -- (3, -3); 

\fill (0, 0) circle (2pt); % Dot at (0, 0)
  \foreach \x in {-2,0,2} {
    \fill (\x, -2) circle (2pt);
  }

\node[above] at (0,0.2) {$\mathbf{1}$};
\node[above] at (0,-1.8) {$h_{-1}$};
\node[above] at (2.2,-1.8) {$e_{-1}$};
\node[above] at (-2.2,-1.8) {$f_{-1}$};    
%%%%%%%%
% Verma picture
% \def\xsh{10};
%   \draw[-] (0+\xsh, 0) -- (-3+\xsh, -3); 
%   \draw[-] (0+\xsh, 0) -- (3+\xsh, -3); 
%   \draw[dashed] (0+\xsh, 0) -- (-5+\xsh, 0); 
  
% \fill (0+\xsh, 0) circle (2pt); % Dot at (0, 0)
%   \foreach \x in {-4,-2,0,2} {
%     \fill (\x+\xsh, -2) circle (2pt);
%   }

%   \foreach \x in {-4,-2} {
%     \fill (\x+\xsh, 0) circle (2pt);
%   }

% \node[above] at (0+\xsh,0.2) {$\mathbf{1}$};
% \node[above] at (0+\xsh,-1.8) {$h_{-1}$};
% \node[above] at (0+\xsh,-2.8) {$e_{-1} f_0$};
% \node[above] at (2.2+\xsh,-1.8) {$e_{-1}$};
% \node[above] at (-2.2+\xsh,-1.8) {$f_{-1}$};

% \node[above] at (-2+\xsh,0.2) {$f_0$};
% \node[above] at (-4+\xsh,0.2) {$f_0^2$};
\end{tikzpicture}
\caption{The vertex algebra $V^k(\mathfrak{sl}_2)$.}
\label{fig:sl2.vac.ver}
\end{figure}
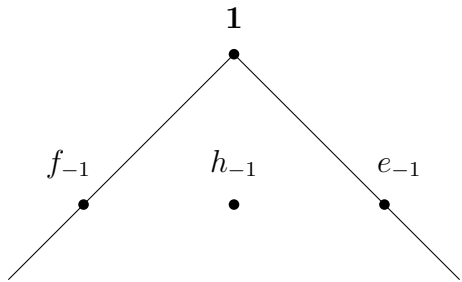

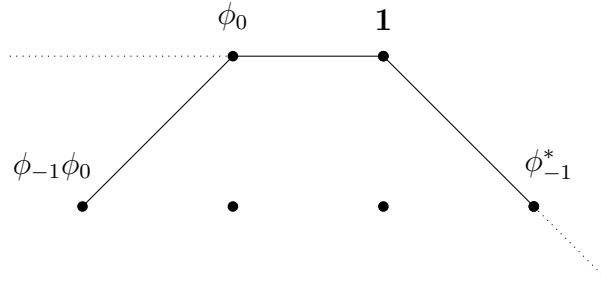
\begin{figure}[H]
\centering
\begin{tikzpicture}
  \draw[-] (0, 0) -- (-2, 0); % Horizontal line going left
    \draw[dotted] (-2, 0) -- (-5, 0); % Horizontal line going left
  \draw[-] (0, 0) -- (2, -2); % Diagonal line going down and to the right
  \draw[dotted] (2, -2) -- (3, -3);
  \draw[-] (-2, 0) -- (-4, -2);
  \fill (0, 0) circle (2pt); % Dot at (0, 0)
  
  \foreach \x in {-2,0} {
    \fill (\x, 0) circle (2pt); % Dots on the horizontal line
  }

\fill (2, -2) circle (2pt);

\node[above] at (0,0.2) {$\mathbf{1}$};
\node[above] at (-2,0.2) {$\phi_0$};
\node[above] at (2.2,-1.8) {$\phi_{-1}^*$};
\node[above] at (-4.4,-1.8) {$\phi_{-1}\phi_0$};

  \foreach \x in {2,0,-2,-4} {
    \fill (\x, -2) circle (2pt);
  }
\end{tikzpicture}
\caption{The rank $1$ charged free fermions vertex algebra $F^{\text{ch}}$.}
\label{fig:fermions}
\end{figure}
We have seen in the previous section that the vertex algebra BRST complex associated with a nilpotent element is given by a tensor product of the form $C = V \otimes F$ where $V = V^k(\g)$ and $F = F^{\text{ch}}(\n)$. In fact the homological $\Z$-grading on $C$ is precisely that inherited from the charge grading on $F$. 

If the differential $d$ in the complex $C$ is compatible with the conformal grading, i.e., if $[L_0, d] = 0$, then we have
\begin{align}\label{eq:Poincare}
\sum_{n \in \Z} \chi_{H^n(C(M))}(q) = \sum_{n \in \Z} \chi_{C^n(M)}(q)
\end{align}
by the Euler-Poincar\'{e} principle (we work with superdimensions here, so components in odd homological degree are odd vector superspaces, and have negative contribution to the superdimension). Later we will see that under favourable circumstances
\[
H^{\neq 0}(C) = 0,
\]
which implies that the character of the BRST reduction can be computed explicitly from knowledge of the character of $V^k(\g)$. Similar remarks apply in the case of nontrivial $V^k(\g)$-modules and their BRST reductions.

It turns out that the zero mode of the most obvious conformal vector on $C$, namely $L \otimes \mathbf{1} + \mathbf{1} \otimes L^{\text{ch}}$, does \emph{not} commute with $d$, but an easy modification fixes the problem: If we set
\begin{align}\label{eq:LC.def}
L^C = \left(L + \frac{1}{2} Th \right) \otimes \mathbf{1} + \mathbf{1} \otimes L^{\text{ch}}
\end{align}
then $[L^C_0, d] = 0$. A similar manoeuvre can be performed in the general case. It transpires that in order to obtain a conformal structure on $C$ compatible with the differential $d = (\Qst + p)_{(0)}$, we should take
\[
L^{\text{ch}} = \sum_{i} \left( \heit(\phi_i) :(T\phi_i) \phi_i^*: + (\heit(\phi_i)-1) :\phi_i (T\phi_i^*): \right),
\]
where $\{\phi_i\}$ is a basis of $\n$ homogeneous with respect to the decomposition \eqref{eq:Dynkin.grading}, and we have written $\heit(\phi) = j$ for $\phi \in \g_j$. This is forced on us by the requirement that $\Delta(p) = 1$. And we should take $L^C$ as in formula \eqref{eq:LC.def}, where now $h \in \g$ is the semisimple element of chosen $\mathfrak{sl}_2$-triple for $f$.

Returning to the case of $\mathfrak{sl}_2$, the conformal vector \eqref{eq:LC.def} gives $e, h, f \in V$ conformal weights $0$, $1$ and $2$ respectively. This is bad news for our computation of the character of $H^0(C)$ via \eqref{eq:Poincare}, because now the character of $V$ has become ill-defined; the infinite set of vectors $e_{-1}^N$ all have the same conformal weight of $0$.
\begin{figure}[H]
\centering
\begin{tikzpicture}
\def\xsh{0};
  \draw[-] (0, 0) -- (5, 0); 
  \draw[-] (0, 0) -- (-2.5, -5); 
  
\fill (0, 0) circle (2pt); % Dot at (0, 0)
  \foreach \x in {0,2,4} {
    \fill (\x, -2) circle (2pt);
  }

  \foreach \x in {2,4} {
    \fill (\x, 0) circle (2pt);
  }

  \foreach \x in {-2,0,2,4} {
    \fill (\x, -4) circle (2pt);
  }

\node[above] at (0,0.2) {$\mathbf{1}$};
\node[above] at (0,-1.8) {$h_{-1}$};
\node[above] at (-2.2,-3.8) {$f_{-1}$};

\node[above] at (2,0.2) {$e_{-1}$};
\node[above] at (4,0.2) {$e_{-1}^2$};
\end{tikzpicture}
\caption{The vertex algebra $V^k(\mathfrak{sl}_2)$ with modified conformal grading.}
\label{fig:sl2.vac.shift}
\end{figure}
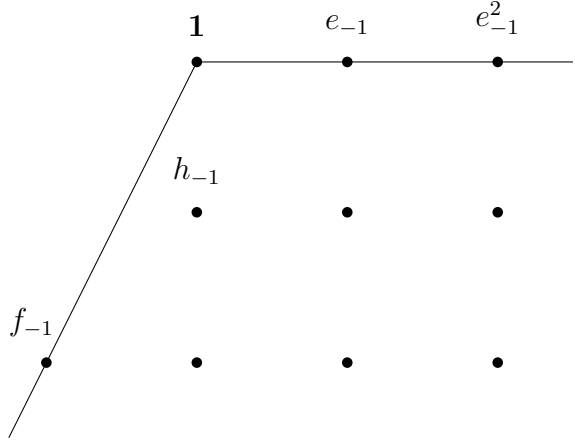
We can obtain a well-defined character by introducing a second gradation, in which $e$, $h$ and $f$ have degrees $-1$, $0$ and $1$, respectively. In other words, gradation by eigenvalues of $-x_0$ where $x = h/2$. The character of $V$ becomes 
\begin{align}\label{eq:V.char.2var}
\chi_V(q) = \tr_V y^{-x_0} q^{(L+Tx)_0} 
&= \prod_{n=1}^\infty \frac{1}{(1- y^{-1} q^{n-1}) (1- q^{n}) (1 - y q^{n+1})}.
\end{align}

Upon passage to homology, the infinite dimensional components of $V$ will be reduced to finite dimension, and $H^0(C)$ will turn out to be conformal with finite dimensional graded pieces (see Theorem \ref{thm:univ.W.str} below). At the level of characters, this fact corresponds to a cancellation between the factor $1-y^{-1}$ in the expression \eqref{eq:F.char.2var} for the character of $F$ (which tends to zero as $y$ tends to $1$), with the factor $(1-y^{-1})^{-1}$ in the expression \eqref{eq:V.char.2var} for the character of $V$ (which tends to infinity as $y$ tends to $1$). Let us then set
\[
\wtil{x} = x \otimes \mathbf{1} + \mathbf{1} \otimes :\phi \phi^*: \in V \otimes F
\]
so that
\begin{align*}
\chi_C(q, y) = \tr_C y^{-\wtil{x}_0} q^{L^C_0} 
&= \prod_{n=1}^\infty \frac{(1-y^{-1} q^{n}) (1-yq^{n-1})}{(1- y^{-1}q^{n+1})(1-q^n)(1-yq^{n-1})} \\
&= (1-yq) \prod_{n=1}^\infty \frac{1}{1-q^n}.
\end{align*}
Finally we can pass to the $y \rightarrow 1$ limit and obtain the graded dimension of the BRST reduction
\[
\chi_{H^0(C)}(q) = \prod_{n=2}^\infty \frac{1}{1-q^n}.
\]
This is precisely the character of the universal Virasoro vertex algebra, and in fact
\[
H^0(C) \cong \vir^{c}
\]
at central charge
\[
c = -\frac{6k^2 + 11k + 4}{k+2}.
\]
(Again, see Theorem \ref{thm:univ.W.str} below.) There is one problem in the story outlined so far, namely that the zero mode $\wtil{x}_0$ does not commute with $d$. It turns out that $\wtil{x}_0$ does commute with $\Qst_{(0)}$, and the required character formula for $H^0(C)$ can be obtained by splitting $d$ as $d = \Qst_{(0)} + p_{(0)}$ and considering the corresponding spectral sequence. Already after passage to the $\Qst_{(0)}$-homology the $L^C_0$ graded pieces become finite dimensional, after which point the limit $y \rightarrow 1$ can be taken legitimately.

Notice that $:\phi \phi^*:$ above is precisely $F^x$, where we are thinking of $\n = \C \phi$ as a $\h$-module.

In general we consider
\[
F^h = \sum_{i} :(\ad(h)\phi_i) \phi_i^*: \quad \text{and} \quad \wtil{h} = h \otimes \mathbf{1} + \mathbf{1} \otimes F^h,
\]
and the arguments go through as above.

The following theorem was proved by Feigin and Frenkel \cite{FF} for $f$ a principal nilpotent element, and in the general case by Kac and Wakimoto \cite{KW04}.
\begin{thm}[{\cite[Theorem 4.1]{KW04}}]\label{thm:univ.W.str}
The BRST homology $H^i(C, d)$ vanishes for $i \neq 0$, and for $i=0$ is a vertex algebra with a set of strong generators parametrised by a basis of $\g^f = \{x \in \g \mid [f, x] = 0\}$. We denote it
\[
W^k(\g, f).
\]
Furthermore, if $\{a_1, \ldots, a_r\}$ is a basis of $\g^f$ homogeneous with respect to the grading $\g = \bigoplus_{j \in \Z} \g_j$ and $a_i \in \g_{-j}$ then the corresponding generator of $W^k(\g, f)$ has conformal weight $j+1$.
\end{thm}

Here is a sketch of the proof, following the approach originated by de Boer and Tjin \cite{dBT} in the case of $f$ principal and $\g = \mathfrak{sl}_n$, and extended to other cases in \cite{FBZ} and \cite{KW04}. One introduces subcomplexes
\[
C_0 = \left< \widetilde{\h + \n_-}, \phi^* \right> \qquad \text{and} \qquad C' = \left< \widetilde{\n}, \phi \right>
\]
of $C = V^k(\g) \otimes F(\n)$, where $\wtil{\n}$ stands for $\{\wtil{x} \mid x \in \n\}$, etc.\footnote{In fact $\n$ is not a $\n_-$-module, so the meaning of $F^{x}$ for $x \in \n_-$ requires explanation. A basis of $\n$ is extended to a basis of $\g$, the action of $\n_-$ is written relative to this basis, and components of the action of $\ad(x)$ in basis vectors outside $\n$ are discarded.} It turns out that $C = C_0 \otimes C'$ as complexes. It is easy to show that $H(C')$ is $\C$ concentrated in degree $0$. So we have
\[
H(C) = H(C_0).
\]
To continue, we  split $d = \Qst_{(0)} + p_{(0)}$, compute cohomology first with respect to $p_{(0)}$, and then with respect to $\Qst_{(0)}$. A further factorisation occurs, with $(C_0, p_{(0)}) = V^\phi(\g^f) \otimes \wtil{C}_0$ as complexes, where the cohomology of $\tilde{C}_0$ is $\C$ concentrated in degree $0$. The spectral sequence collapses on the next page, and we obtain
\[
H^0(C) = V^\phi(\g^f)
\]
as (graded) vector spaces, as required.

Determining the vertex algebra structure of $W^k(\g, f)$ explicitly requires lifting the elements of $V(\g^f)$ to $d$-cocycles, which turns out to be very complicated in general. Nevertheless, we have the character of $W^k(\g, f)$; it is
\[
\prod_{j \in \Z_+} \prod_{n=1}^\infty \frac{1}{1-q^{n + \dim(\g^f_{-j}) + 1}}
\]

The algebra of functions on the Slodowy slice $S_f$ is the BRST reduction of the Poisson algebra $S(\g)$ by a Hamiltonian action. Replacing $S(\g)$ by its obvious quantisation $U(\g)$ and carrying out the BRST recipe, yields a class of associative algebras $U(\g, f)$ known as finite $W$-algebras. One might ask for the relation between $W^k(\g, f)$ and $U(\g, f)$. The answer turns out to be as simple as one could hope for:
\begin{thm}[{\cite{DSK06}}]
The Zhu algebra of $W^k(\g, f)$ is the finite $W$-algebra $U(\g, f)$.
\end{thm}
The case of principal nilpotent element $f$ had been proved earlier by Arakawa \cite{Ar.rep.I}, in this case the finite $W$-algebra coincides with the centre $Z(\g)$ of the associative algebra $U(\g)$.

The argument, in a nutshell, is this. The complex $C = V^k(\g) \otimes \Fch(\n)$ is itself a vertex superalgebra, and so has a Zhu algebra. This associative algebra inherits a differential, and there is a natural map
\[
\zhu(W^k(\g, f)) \rightarrow H^0(\zhu(C^\bullet), d).
\]
One directly verifies that the Zhu algebra complex is nothing other than the finite BRST complex which computes $U(\g, f)$ by the definition of the latter.

There is a complication, however. The preceding statements are true of the vertex algebra $V^k(\g) \otimes \Fch(\n)$ in which $V^k(\g)$ is equipped with its standard conformal vector.  But we were forced to modify this vector
\[
L \rightarrow L + Tx
\]
to ensure compatibility with $d$. This change alters the conformal weights and conformal gradings of the quantum fields in a systematic way, specifically
\[
a_n \mapsto a_{n-j} \quad \text{for $a \in \g_j$},
\]
and thus alters the Zhu algebra. The effect of this modification is ultimately that the Zhu algebra is still isomorphic to $U(\g, f)$, but the bijection between irreducible positive energy $W^k(\g, f)$-modules and irreducible $U(\g, f)$-modules is more convoluted than one might expect.

We continue in the principal case for simplicity. Irreducible $Z(\g)$-modules are one dimensional and are known as central characters. For example any Verma $\g$-module $M(\lambda) = U(\g) v_{\la}$ induces a central character
\[
\gamma_\lambda : Z(\g) \rightarrow \en(\C)
\]
defined by $z v_\lambda = \gamma_\lambda(z) v_\lambda$ for all $z \in Z(\g)$. Chasing through the identifications sketched above, one discovers that the $Z(\g)$-module corresponding to the principal affine $W$-algebra itself is
\[
W^k(\g, f_{\text{prin}})_{\low} \cong \gamma_{-(k+h^\vee)x}.% = \gamma_{-(k+h^\vee)\rho^\vee}.
\]

We now describe a variant of the BRST reduction functor, which has superior formal properties for applications in representation theory. We build a charged free fermion vertex superalgebra based on $\n_-$ now, that is $F_- = F^{\text{ch}}(\n_-)$. We equip $F_-$ with a conformal structure such that $\Delta(\n_-) = 0$ and $\Delta(\n_-^*) = 1$.

Let $M$ be a $V^k(\g)$-module. We define $H_-(M)$ to be the cohomology of the complex
\[
C_-(M) = M \otimes F_-, \quad d = (Q^- + p_-)_{(0)}
\]
where
\begin{align}\label{eq:C.minus.def}
\begin{split}
Q^- &= \sum_{i \in \D_-} \alpha(\phi_i) \otimes \phi_i^* - \mathbf{1} \otimes \frac{1}{2} \sum_{i, j \in \D_-} : [\phi_i, \phi_j] \, \phi_i^* \phi_j^*: \\
p_- &= \sum_i (e, \phi_i) \mathbf{1} \otimes \phi_i^*.
\end{split}
\end{align}

The vector space $H_-(M)$ carries a natural $W^k(\g, f)$-module structure. The explicit description of this action involves $w_\circ$ the longest element in the Weyl group of $\g$ and a few other ingredients. We refer the reader to {\cite[p. 291]{Ar.rep.I}} for details.

Arakawa proved the following vanishing theorem for $H_-(-)$. Recall the decomposition $\g = \bigoplus_{j \in \Z} \g_j$ associated with $f$. Recall also that the category of $V^k(\g)$-modules is equivalent to the category of smooth $\what{\g}$-modules of level $k$. Let $\OO_k$ denote category of $\what{\g}$-modules $M$ of level $k$, admitting a weight space decomposition $M = \bigoplus_{\mu \in \what\h^*} M_\mu$ with finite dimensional weight spaces, such that $\supp(M) \subset \cup_{i} (\mu_i - Q_+)$, where $i$ runs over some finite set. Further let $\OO_k^{(0)}$ be the full subcategory of $\OO_k$ consisting of $\what\g$-modules integrable as modules over the subalgebra $\g_0 \subset \g \subset \what\g$. If $f$ is principal then $\OO_k^{(0)} = \OO_k$ because $\g_0 = \h$.
\begin{thm}[{\cite[Theorem 7.6.1]{Ar.rep.I}} {\cite[Theorem 5.4.1]{Ar.rep.II}}]{\ }

Let $\g$ be a simple Lie algebra, $f \in \g$ nilpotent, and other constructions as described above.
\begin{itemize}
\item If $f$ is principal, and $M \in \OO_k$, then $H^{i}_-(M) = 0$ for all $i \neq 0$,

\item If $f$ is even, and $M \in \OO_k^{(0)}$, then $H^i_-(M) = 0$ for all $i\neq 0$,

\item More generally if $f$ possesses a good even grading, and $M \in \OO_k^{(0)}$, then $H^i_-(M) = 0$ for all $i\neq 0$ (though the $H^i_-(M)$ are not $W^k(\g, f)$-modules, rather they are Ramond twisted $W^k(\g, f)$-modules).
\end{itemize}
\end{thm}

Next comes the following theorem, which establishes the structure of $H_-(-)$ reduction of Verma $\what\g$-modules. To state the theorem, we need to elaborate a little on the relation between irreducible positive energy $V$-modules and irreducible $\zhu(V)$-modules. It is possible to construct a Lie algebra $\Lie(V)$ canonically from $V$, whose elements are symbols $a_n$ as $a$ ranges over $V$ and $n$ ranges over $\Z$, with the property that any $V$-module becomes a $\Lie(V)$-module in the obvious way. This Lie algebra is $\Z$-graded and, it turns out, there is a homomorphism $\Lie(V)_0 \rightarrow \zhu(V)$, where $\zhu(V)$ is equipped with the commutator bracket. Starting from a $\zhu(V)$-module $N$, which we can convert into a $\Lie(V)_0$-module, we may impose an action of $\Lie(V)_{>0}$ in $N$ by $0$, and then induce to obtain a $U(\Lie(V))$-module
\[
\wtil{M} = U(\Lie(V)) \otimes_{U(\Lie(V))_{\geq 0}} N.
\]
A quotient of $\wtil{M}$ by a subspace $(B)\cdot\wtil{M}$ (which corresponds to the imposition of Borcherds identity), results in a $V$-module
\[
\mathbb{M}(N) = \wtil{M} / (B)\cdot\wtil{M}
\]
which is positive energy with $\mathbb{M}(N)_{\low}$. The irreducible module $\LLL(N)$, already introduced in Section \ref{sec:lec2}, is the quotient of $\mathbb{M}(N)$ by its maximal proper $V$-submodule. The relationship between $\mathbb{M}(N)$ and $\LLL(N)$ is reminiscent of the relationship between Weyl $\what\g$-modules and irreducible highest weight $\what\g$-modules.
\begin{thm}[{\cite[Theorem 7.5.1]{Ar.rep.I}}]
The reduction of Verma modules
\[
H_{-}^0(M_k(\lambda)) \cong \mathbb{M}(\gamma_{\lambda^*})
\]
Here $\lambda^* = -w_\circ(\lambda)$ is the dual weight.
\end{thm}

To prove this theorem, one first checks that the central character of the lowest piece is as stated. Note that $\lambda$ is ``twisted'' to $\lambda^*$, this is due to the structure of the action of $W^k(\g, f)$ on $H_-(M)$ as commented after equation \eqref{eq:C.minus.def}. By an analysis similar to the proof of Theorem \ref{thm:univ.W.str}, it is shown that $H_-^0(M(\lambda))$ has a PBW basis. This implies that the canonical map
\[
\mathbb{M}(\gamma_{\lambda^*}) \rightarrow H_{-}^0(M_k(\lambda))
\]
must be an isomorphism.

Finally, using the previous theorem and exactness of $H_-^0(-)$, Arakawa proves the following theorem.
\begin{thm}[{\cite[Theorem 7.6.3]{Ar.rep.I}}]
For $\lambda$ antidominant, the reduction of $L_k(\lambda)$ is
\[
H_{-}^0(L_k(\lambda)) \cong \mathbb{L}(\gamma_{\lambda^*})
\]
Here $\lambda^* = -w_0 \lambda$ is the dual weight.
\end{thm}
This result really is key since with it, knowledge of the characters of $V^k(\g)$-modules can be used to compute the characters of the corresponding $W^k(\g, f)$-modules.

We sketch the proof. By exactness of $H_-^0(-)$ we have a surjection
\[
\mathbb{M}(\gamma_{\lambda^*}) = H^0_-(M_k(\lambda)) \rightarrow H_-^0(L_k(\lambda)),
\]
and an injection
\[
H^0_-(L_k(\lambda)) \rightarrow H^0_-(D(M_k(\lambda))) = D(\mathbb{M}(\gamma_\lambda)).
\]
A standard argument (see \cite[Theorem 3.3]{Humphreys} for the prototype case of dual Verma modules over finite dimensional simple Lie algebras) says that image of $\mathbb{M}(\gamma_{\lambda^*})$ in $D(\mathbb{M}(\gamma_\lambda))$ is the unique simple highest weight submodule of the latter, isomorphic to $\LLL(\gamma_{\la^*})$.

\section{Exceptional $W$-algebras}

The universal affine vertex algebra $V^k(\g)$ is simple for generic values of the level $k$, but at special levels $V^k(\g)$ has a nontrivial simple quotient $L_k(\g)$ whose representation theory is quite interesting. In Theorem \ref{thm:aff.k.Zplus} for instance we saw that $L_k(\g)$ has finitely many irreducible positive energy modules when $k \in \Z_+$. Furthermore in this case $L_k(\g)$ is an example of a \emph{rational} vertex algebra, which loosely means that all of its modules are completely reducible into direct sums of irreducible modules, and that it has finitely many irreducible modules. For admissible level $k$ the vertex algebra $L_k(\g)$ is not rational, but we have seen in Theorem \ref{thm:A.rat.O} that some parallels with the positive integral case persist.

On the other hand we have examined the BRST procedure, and obtained a class of vertex algebras $W^k(\g, f)$, the universal affine $W$-algebras. Once again these vertex algebras are simple for generic level $k$, but have interesting simple quotients at special levels. A basic question is: for which combinations of $\g$, $f$ and $k$ is the simple quotient $W_k(\g, f)$ of $W^k(\g, f)$ a rational vertex algebra?

A natural class of candidates is given by the \emph{exceptional $W$-algebras}, which we define below. In this section we describe some results from our article \cite{AVE-EMS} on the structure and representation theory of exceptional $W$-algebras, including rationality in many cases.

In Theorem \ref{thm:aff.k.Zplus} we saw the classification of irreducible positive energy $L_k(\g)$-modules, for $k \in \Z_+$. By the correspondence with irreducible modules over the Zhu algebra, it follows that
\[
\zhu(L_k(\g)) \cong \bigoplus_{\lambda \in P_+^k} \en(L(\la)).
\]
In \cite{AVE-EMS} we generalise this result from positive integral to other admissible levels. To state the result, we first fix some notation. Let $\g$ be a finite dimensional simple Lie algebra, and $U = U(\g)$. For a $U$-module $M$ we consider the two-sided ideal
\[
\ann_U(M) = \{u \in U \mid u M = 0\},
\]
and for $\la \in \h^*$ we write $J_\la$ for $\ann_U(L(\la))$. Recall that a weight $\la$ is said to be dominant if $\left<\la+\rho, \al^\vee\right> \notin \Z_{\leq 0}$ for all $\al \in \D_+$. A result of Joseph asserts that if $\la$ and $w \circ \la$ are dominant then $J_{w \circ \la} = J_\la$. In particular admissible weights are dominant, so if we write $\Pr_k$ for the set of admissible weights of admissible level $k$, then $J_\la$ depends on $\la \in \Pr_k$ only through its image in
\[
[{\Pr}_k] = {\Pr}_k / \sim,
\]
where $\sim$ denotes the equivalence relation on $\h^*$ induced by the $W \circ (-)$ action.
\begin{thm}[{\cite{AVE-EMS}}]
Let $k$ be an admissible level for $\g$, then 
\[
\zhu(L_k(\g)) \cong \prod_{\la \in [{\Pr}_k]} U / J_\la.
\]
\end{thm}

The Zhu algebra of $L_k(\g)$ (for any level $k$) is isomorphic to the quotient $U / I_k$ for some ideal $I_k$ which depends on $k$. We consider the tensor product of $\zhu(L_k(\g))$, over the centre $Z = Z(U)$, with a central character. The basic observation in the proof of the theorem is that
\begin{align*}
\zhu(L_k(\g)) \otimes_Z \C \gamma_\lambda = U / J
\end{align*}
for some two sided ideal $J$, and since $M(\lambda) / J M(\lambda)$ would be a highest weight $\zhu(L_k(\g))$-module, it can only be nonzero if $\la \in {\Pr}_k$ and furthermore we must have $M(\lambda) / J M(\lambda) = L(\la)$, and hence $J = J_\lambda$.

Recall that $\zhu(W^k(\g), f) \cong U(\g, f)$. Since $L_k(\g)$ is a quotient of $V^k(\g)$, the Zhu algebra of $H^0(L_k(\g))$ is a certain quotient of $U(\g, f)$, which we would now like to determine via BRST reduction, using the description obtained above of the Zhu algebra of $L_k(\g)$. Here and below let us write
\[
W = H^0(L_k(\g)).
\]
\begin{thm}[{\cite{AVE-EMS}}]
Let $k$ be an admissible level for $\g$, then 
\[
\zhu(W) \cong \prod_{\la \in [{\Pr}_k]} H^0(U / J_\la).
\]
\end{thm}
Note that it can easily happen that $W = 0$. Indeed if $k \in \Z_+$ then $L_k(\g)$ is $G$-integrable and so $H^0(L_k(\g)) = 0$ can already be deduced from consideration of characters as in Section \ref{sec:str-mod-char}. If $W \neq 0$ then certainly there exists a surjection from $W$ to $W_k(\g, f)$ the unique simple quotient.

The universal enveloping algebra $U = U(\g)$ is naturally filtered by ``word length'' and the associated graded algebra $\gr(U)$ is isomorphic to the symmetric algebra $S(\g)$ (in fact this is one formulation of the PBW theorem). Let $\g$ be finite dimensional and simple, so that $S(\g)$ is identified with the algebra of  polynomials on $\g \cong \g^*$. An ideal $I \subset U$ induces an ideal $\gr(I) \subset S(\g)$ and the corresponding zero set is known as the associated variety $\Var(I) \subset \g$. We may also write, more specifically, $\Var(M)$ for $\Var(\ann_U(M))$.

Let $k$ be an admissible level for $\g$, then a theorem of Arakawa {\cite[Theorem 5.14]{Arakawa.2015}} asserts that $\Var(I_k) \subset \g$ is in fact the closure of a nilpotent orbit, which depends only on the denominator, so we call it $\mbo_q$.

Now a result of Losev \cite{Losev} asserts that
\begin{align*}
H^0(U / J_\lambda) =
\begin{dcases}
0 & \Var(J_\lambda) \subsetneq \ov{G \cdot f} \\
\text{finite dimensional} & \Var(J_\lambda) = \ov{G \cdot f}
\end{dcases}
\end{align*}
Of course, for all $\lambda \in {\Pr}_k$ we have
\[
\Var(J_\lambda) \subset \Var(I_k).
\]
%Call $Pr_\circ^k$ those $\lambda$ satisfying this condition.
It therefore makes sense to consider the restriction on our choice of $f$:
\begin{hypoth}\label{hypoth:1}
The nilpotent element $f$ is chosen in the open part of $\Var(I_k)$, namely in $\mbo_q$.
\end{hypoth}
In fact let us elevate this hypothesis to a definition
\begin{defn}
The $W$-algebra $H^0(L_k(\g)), f)$ will be said to be exceptional if $f \in \mbo_q$.
\end{defn}
(Actually this definition is similar but slightly different to the notion of exceptional introduced in \cite{KW09}). One of the key results of \cite{Arakawa.2015} is that exceptional $W$-algebras are $C_2$-cofinite. In general a vertex algebra $V$ is said to be $C_2$-cofinite (or ``lisse'') if $V / V_{(-2)}V$ is finite dimensional. The importance of this condition, which was introduced by Zhu in \cite{zhu96}, is certainly not immediately clear, but it turns out to imply several other desirable finiteness and convergence properties for objects associated with $V$. In particular, although $C_2$-cofiniteness is known not to imply rationality, it suggests it. 

Anyway, under Hypothesis \ref{hypoth:1}, we have $\Var(I_k) = \ov{G \cdot f}$, and so $\Var(J_\lambda) \subset \ov{G \cdot f}$ for all $\lambda \in {\Pr}_k$. Let us call an admissible weight $\lambda$ ``replete'' if we have equality $\Var(J_\lambda) = \ov{G \cdot f}$. Intuitively the replete weights $\la$ are those for which we expect a nonzero $W$-module obtained as Hamiltonian reduction of $L(\la)$.

\begin{cor}[{\cite{AVE-EMS}}]
The Zhu algebra of an exceptional $W$-algebra is finite dimensional and semisimple.
\end{cor}

\begin{remark}
McRae has recently proved that a $C_2$-cofinite vertex algebra $V$ for which $\zhu(V)$ is finite dimensional and semisimple (and for which a few other more technical hypotheses are satisfied) is rational \cite{McRae-preprint}. It follows that exceptional $W$-algebras are rational.
\end{remark}
In order to go further and explicitly describe the irreducible $W$-modules using Hamiltonian reduction, we have found it convenient to impose some technical conditions. (For some progress in cases outside the scope of these conditions, see \cite{Fasquel}.)

Let us call an admissible weight $\lambda$ ``conservative'' for $f$ if
\begin{itemize}
\item $\lambda$ is dominant integral with respect to $\g_0$ (i.e. $L(\lambda) \in \OO_k^{(0)}$), and

\item $\lambda$ is antidominant with respect to all roots of $\D_{>0}$ (i.e., $\left<\lambda+\rho, \al^\vee\right> \notin \Z_{>0}$ for all $\al \in \D_{>0}$).
\end{itemize}
Then we have
\begin{prop}[{\cite{AVE-EMS}}]\label{prop:technical.goldilocks}
If $\lambda$ is conservative for $f$, then the associative algebra $H^0(U / J_\lambda)$ has a unique irreducible module $E_\lambda$, and
\[
H^0_{-}(L_k(\lambda)) \cong \LLL(E_\la).
\]
\end{prop}

We sketch the proof. Results of \cite{Ar.rep.II}, specifically Theorem 5.5.4 guarantee that on the category $\OO_k^{(0)}$ of $\g_0$-integrable $\what\g$-modules, the reduction functor $H_{-}^0(-)$ is exact and sends conservative modules to nonzero modules. Furthermore that $H_-(L_k(\lambda))$ is almost irreducible, meaning that the union of all its proper submodules has trivial intersection with $H_-(L_k(\lambda))_{\low}$. This produces a nonzero $H^0(U / J_\lambda)$-module, which means in particular that this algebra must be nonzero, and so we deduce $\Var(J_\lambda) = \overline{G \cdot f}$. On the other hand, Losev characterised the set of finite dimensional irreducible $H^0(U / J_\lambda)$-modules: they are finite in number and there is a certain group $C(f)$ which acts transitively on this set. Integrability of $L(\lambda)$ under $G_0$ (the algebraic group corresponding to $\g_0$) assures us that $C(f)$ acts trivially, and so $H^0_{-}(L_k(\lambda))_{\low}$ must be the unique irreducible module, hence the isomorphism of the theorem statement.

Proposition \ref{prop:technical.goldilocks} gives us good control of the Hamiltonian reduction of $L_k(\lambda)$ for conservative $\lambda$. We therefore introduce the following:
\begin{hypoth}
The $W \circ (-)$-class of every replete weight has a conservative representative.
\end{hypoth}
Under this hypothesis, the representation theory of $W$ is sufficiently well controlled by that of $L_k(\g)$, that we may prove the following rationality theorem.
\begin{thm}[{\cite{AVE-EMS}}]
Under the hypothesis above, the exceptional $W$-algebra is rational.
\end{thm}

We sketch the proof. It is not difficult to argue that the existence of a nontrivial extension
\[
0 \rightarrow \LLL(E_{\lambda}) \rightarrow M \rightarrow \LLL(E_{\lambda'}) \rightarrow 0,
\]
implies $[M_k(\lambda') : L_k(\lambda)] > 0$, which in turn implies linkage of $\lambda$ and $\lambda'$ under the admissible level affine Weyl group $W \ltimes p Q^\vee$. But dominance of both $\la$ and $\la'$ then implies that they are actually equal. On the other hand self-extensions of $\LLL(E_\lambda)$ would be detected by the Zhu algebra, and are thus ruled out since it is semisimple.

\section{Characters of exceptional $W$-algebras}\label{sec:excep.chars}

We have seen that application of the Weyl-Kac character formula to integrable highest weight $\what\g$-modules yields character formulas in terms of classical theta functions. This phenomenon has its origin in the structure of the affine Weyl group $\what{W} = W \ltimes t_{Q^\vee}$ as a semidirect product of a finite group with a lattice, together with the structure of the action of the lattice $t_{Q^\vee}$ on weights, see formula \eqref{eq:transl.k}. We have also seen that the same goes through for admissible highest weight modules $L(\lambda)$, the group $\what{W}$ being replaced by its subgroup $\what{W}(\lambda)$ which has substantially the same structure as $\what{W}$.

Let $Q$ be an integral lattice of rank $\ell$, and $\h = Q \otimes_\Z \C$. By a theta function of $Q$ we mean a function of the general form
\begin{align}\label{eq:lattice.theta}
\theta_{\mu+Q}(\tau, x) = \sum_{\alpha \in \mu + Q} e^{2\pi i (\alpha, x)} e^{\pi i \tau (\alpha, \alpha)}, 
\end{align}
defined on the domain $(\tau, x) \in \HH \times \h$ where $\HH$ is the upper half complex plane. A key feature of theta functions is their modular invariance. In particular, for each $\mu + Q \in Q^\vee / Q$ one has
\[
\theta_{\mu+Q}(\tau, z)(-1/\tau, x/\tau) = \tau^{\ell/2} e^{\pi i |x|^2 / \tau} \sum_{\mu' \in Q^\vee / Q} S_{\mu, \mu'} \theta_{\mu'+Q}(\tau, x),
\]
where
\[
S_{\mu, \mu'} = \frac{(-i)^{\ell/2}}{|Q^\vee / Q|^{1/2}} e^{-2\pi i (\mu, \mu')}.
\]

Evaluation of the Kac-Wakimoto character formula \eqref{eq:KW.char} on $2\pi i \tau d + x \in \h^*$ yields an expression
\[
\chi_{L(\lambda)}(\tau, x) = B_\lambda(\tau, x) / R(\tau, x),
\]
where
\begin{align}\label{eq:R.theta}
R(\tau, x) = \prod_{\alpha \in \D_+} \Theta(\tau, \alpha(x)), \qquad \Theta(\tau, z) = \prod_{n=1}^{\infty} (1-e^{-2\pi i z}q^{n-1})(1-q^n)(1-e^{2\pi i z}q^n)
\end{align}
is the Weyl denominator, and $B_\lambda(\tau, x)$ is a sum of theta functions \eqref{eq:lattice.theta} for a rescaled copy of the root lattice $Q$ of $\g$, the sum being parametrised by the finite Weyl group $W$. 

We now examine two cases: the case of principal nilpotent element $f$, and the case of subregular nilpotent element $f$.

Let $k = -h^\vee + p/q$ be an admissible weight for $\g$. Recall that an exceptional $W$-algebra $H^0(L_k(\g))$ is one for which the nilpotent element $f$ lies in the open part of $\Var(I_k) \subset \g$, which is to say, $f \in \mbo_q$. It happens that $\mbo_q$ equals the principal nilpotent orbit for all sufficiently large $q$, namely $q \geq h^\vee$. So let us assume that $q \geq h^\vee$ throughout this section.

We now describe the admissible weights of level $k$. In fact the first part of this discussion is valid for any admissible level $k$, without the assumption $q \geq h^\vee$. It is not hard to see that the admissible weights $\la$ of level $k$ whose integral coroot system is $\Pi^\vee(\la) = \Pi^\vee_q$ consists of
\[
\lambda = k\Lambda_0 + \nu, \quad \text{for $\nu \in P_+^{p, \text{reg}}$}.
\]
Here $P_+^{p, \text{reg}}$ means the set of weights
\[
\{\nu \in P \mid \text{$\left< \nu, \al_i^\vee \right> > 0$ for $i=1, \ldots, \ell$ and $\left< \nu, \theta^\vee \right> < p$} \}
\]
\begin{figure}[H]
\centering
%\begin{tikzpicture}[scale=0.4, decoration={markings,mark=at position 1 with {\arrow[scale=1.5,black]{latex}};}]
\begin{tikzpicture}[scale = 0.6]
\draw[fill=black] (-3cm,-1.73cm) circle (.05cm);
\draw[fill=black] (-1cm,-1.73cm) circle (.05cm);
\draw[fill=black] (1cm,-1.73cm) circle (.05cm);
\draw[fill=black] (3cm,-1.73cm) circle (.05cm);
\draw[fill=black] (5cm,-1.73cm) circle (.05cm);
\draw[fill=black] (7cm,-1.73cm) circle (.05cm);
\draw[fill=black] (9cm,-1.73cm) circle (.05cm);
\draw[fill=black] (-2cm,0cm) circle (.05cm);
\draw[fill=black] (0cm,0cm) circle (.05cm);
\draw[fill=black] (2cm,0cm) circle (.05cm);
\draw[fill=black] (4cm,0cm) circle (.05cm);
\draw[fill=black] (6cm,0cm) circle (.05cm);
\draw[fill=black] (8cm,0cm) circle (.05cm);
\draw[fill=black] (-1cm,1.73cm) circle (.05cm);
\draw[fill=black] (1cm,1.73cm) circle (.05cm);
\draw[fill=black] (3cm,1.73cm) circle (.2cm);%%
\draw[fill=black] (5cm,1.73cm) circle (.2cm);%%
\draw[fill=black] (7cm,1.73cm) circle (.05cm);
\draw[fill=black] (0cm,3.46cm) circle (.05cm);
\draw[fill=black] (2cm,3.46cm) circle (.05cm);
\draw[fill=black] (4cm,3.46cm) circle (.2cm);%%
\draw[fill=black] (6cm,3.46cm) circle (.05cm);
\draw[fill=black] (1cm,5.20cm) circle (.05cm);
\draw[fill=black] (3cm,5.20cm) circle (.05cm);
\draw[fill=black] (5cm,5.20cm) circle (.05cm);
\draw[fill=black] (2cm,6.92cm) circle (.05cm);
\draw[fill=black] (4cm,6.93cm) circle (.05cm);
\draw[dashed] (60: 0 mm) -- (60: 90 mm);
\draw[dashed] (0: 0 mm) -- (0: 90 mm);
\draw[dashed] (0: 80 mm) -- (60: 80 mm);
\def\x{0.98};
\draw[-{Latex[length=2mm,width=1mm]}] (0,0) -- (3cm*\x,-1.73cm*\x);
\draw[-{Latex[length=2mm,width=1mm]}] (0mm,0mm) -- (0mm,3.46cm*\x);
%%
%% LABELS
%%
\node[] at (-0.8, 3.8) {$\alpha_1$};
\node[] at (3.8, -1.5) {$\alpha_2$};
\node[] at (0.9, 2.6) {$\varpi_1$};
\node[] at (2.5, -0.7) {$\varpi_2$};
\end{tikzpicture}
  \caption{The set $P_+^{p, \text{reg}}$ of regular weights of level $p=4$ for $\g = \mathfrak{sl}_3$}
  \label{fig:reg}
\end{figure}
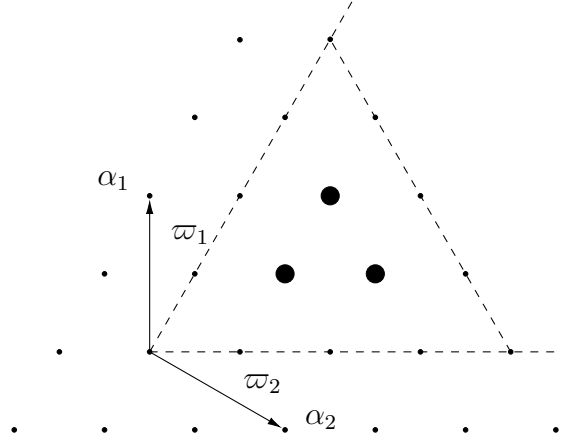

The remaining admissible weights $\la$ are obtained by classifying the possible sets of simple coroots $\Pi^\vee(\lambda)$, obtained as transformations of $\Pi^\vee_q$
\begin{align*}
\Pi^\vee(\lambda) &= y t_{-\eta} (\Pi^\vee_q) \\
\lambda &= yt_{-\eta}(\nu) - \rho,
\end{align*}
where $y \in W$ the finite Weyl group and $\eta \in P$ the weight lattice. Roughly speaking, the transformation $t_{-\eta}$, applied to $\Pi^\vee_q$ lifts the simple coroots $\alpha^\vee_i$, $i=1, \ldots, \ell$ out of $\D^\vee_+$ into $\widehat{\D}^{\vee, \text{re}}_+$, sending $\al^\vee_i \mapsto \al^\vee_i + m_i K$ with the coefficients $m_i \in \Z_{\geq 0}$ encoded into $\eta \in P$. Then $y$ sends $\al^\vee_i + m_i K \mapsto y(\al^\vee_i) + m_i K$. Since the simple coroot $-\theta^\vee + q K$ must remain in $\widehat{\D}^{\vee, \text{re}}_+$ at the end of this process, it is required that $\eta \in P_+^{q}$. In fact one perceives that this condition is not sufficient, and in general the parametrisation of valid pairs $(y, \eta)$ is quite complicated.

We are analysing the case of $f$ principal nilpotent. And in this case a fortuitous simplification takes place. For such $f$, an admissible weight $\lambda$ is replete if it is antidominant, and any replete weight is automatically conservative, because $\g_0 = \h$ and integrability with respect to this subalgebra is no  condition on $\la$. We are thus interested in parametrising those $\la \in {\Pr}_k$ which pair non-integrally with each finite simple coroot $\al^\vee_i$, $i =1, \ldots, \ell$. This condition on $\lambda$ translates directly into the condition $\eta \in P_+^{q, \text{reg}}$. Now since $t_{-\eta}(\D_+^\vee) \subset \widehat{\D}^{\vee, \text{re}}_+ \backslash \D_+$ it follows that $y t_{-\eta} (\Pi^\vee_q) \subset \widehat{\D}^{\vee, \text{re}}_+$ for any $y \in W$.

We thus obtain (following Kac and Wakimoto {\cite[Lemma 4.2]{KW90}}) a pleasingly symmetrical parametrisation of the relevant admissible weights:
\[
(y, \nu, \eta) \in W \times P_+^{p, \text{reg}} \times P_+^{q, \text{reg}}.
\]
The BRST reduction sends $L(\lambda)$ and $L(y \circ \lambda)$ to the same irreducible $W$-algebra module, so the action of the finite Weyl group $W$ is removed from the parametrisation.

The principal nilpotent $f$ can be chosen so that in the decomposition \eqref{eq:Dynkin.grading}, the subalgebra $\n = \g_{>0}$ coincides with the subalgebra $\n_+$ of the usual triangular decomposition $\g = \n_- \oplus \h \oplus \n_+$. A little calculation shows that the character of $\Fch(\n)$ is ``two thirds'' of the Weyl denominator $R$ of \eqref{eq:R.theta}. More precisely
\[
\chi_{\Fch(\n)}(\tau, x) = \prod_{\alpha \in \D_+} \Theta(\tau, \alpha(x)),
\]
where
\[
\Theta(\tau, z) = \prod_{n=1}^{\infty} (1-e^{-2\pi i z}q^{n-1})(1-e^{2\pi i z}q^n).
\]
At the level of characters the BRST procedure converts the character $\chi_{L(\la)}(\tau, x) = B_\la(\tau, x) / R(\tau, x)$ into
\[
\chi_{H(L(\la))}(\tau) = \lim_{x \rightarrow 0} \frac{B_\la(\tau, x)}{\prod_{n=1}^{\infty} (1-q^n)^\ell}.
\]
It is thus possible to explicitly determine the modular transformations of the characters
\[
\chi_{H(L(\la))}(-1/\tau) = \sum_{\la'} S_{\la, \la'} \chi_{H(L(\la'))}(\tau).
\]
This was done by Frenkel, Kac and Wakimoto \cite{FKW}, uncovering a very pleasant factorisation of the $S$-matrix
\begin{align}\label{eq:prin.S.mat}
S_{\la, \la'} = S_{(\nu, \eta), (\nu', \eta')} = \sum_{y \in W} \varepsilon(y) e^{-2\pi i \frac{p}{q} (y(\eta), \eta')} \sum_{w \in W} \varepsilon(w) e^{-2\pi i \frac{q}{p} (w(\nu), \nu')},
\end{align}
related to the Feigin-Frenkel isomorphism \cite{FF}.

We now pass to the case of subregular nilpotent element $f$, and results of our article \cite{AVE-EMS}. The subregular nilpotent orbit in the finite dimensional simple Lie algebra $\g$ is uniquely characterised by the property $\dim(\g_0) = \dim(\h) + 2$ in the decomposition \eqref{eq:Dynkin.grading}. From now on we assume $\g$ to be simply laced, then it turns out that the subregular nilpotent $f$ can be chosen so that $\g_0 = \h + \g_{\al_*} + \g_{-\al_*}$ for a certain root $\al_*$. For types $D$ and $E$ this root $\alpha_*$ corresponds to the trivalent node in the Dynkin diagram. For type $A_n$ the subregular nilpotent $f$ is even only if $n$ is odd, in which case $\alpha_*$ corresponds to the middle node in the Dynkin diagram. The case $n$ even can be treated by introducing a ``good even grading'' on $\g$ adapted to $f$, different than the Dynkin grading. We have avoided discussing this trick in these notes to keep things as simple as possible.

Let $k = -h^\vee + p/q$ be admissible. The denominator(s) $q$ for which $\mbo_q$ equals the subregular nilpotent orbit were determined in \cite{Arakawa.2015}. For example for $D_n$ one has $q = 2n-4$ and $q = 2n-3$, and for $\g = E_8$ one has $q = 24, 25, \ldots, 29$. We need to identify the replete weights and conservative weights $\la \in {\Pr}_k$. Let us return to the discussion around $\Pi^\vee(\la) = y t_{-\eta}(\Pi^\vee_q)$. In the principal case it was necessary that $t_{-\eta}$ shifted all the simple coroots $\al_i^\vee$ out of the finite root system $\D$, so that $\la$ would not pair integrally with any finite root. In the present case, $\la$ is replete if $t_{-\eta}(\al_i^\vee) = \al_i^\vee$ for exactly one $i=1,\ldots, \ell$ or else $t_{-\eta}(-\theta_i^\vee + qK) = -\theta_i^\vee$. And $\la$ is conservative if furthermore $\left<\lambda, \alpha_*\right> \in \Z_{\geq 0}$.

Interestingly, the repleteness condition says precisely that $\eta$, instead of being regular as it was in the principal case, must lie on exactly one of the walls of the level $q$ affine Weyl chamber:
\begin{figure}[H]
\centering
  \begin{tikzpicture}[scale=0.6]
\draw[fill=black] (-3cm,-1.73cm) circle (.05cm);
\draw[fill=black] (-1cm,-1.73cm) circle (.05cm);
\draw[fill=black] (1cm,-1.73cm) circle (.05cm);
\draw[fill=black] (3cm,-1.73cm) circle (.05cm);
\draw[fill=black] (5cm,-1.73cm) circle (.05cm);
\draw[fill=black] (7cm,-1.73cm) circle (.05cm);
\draw[fill=black] (9cm,-1.73cm) circle (.05cm);
\draw[fill=black] (-2cm,0cm) circle (.05cm);
\draw[fill=black] (0cm,0cm) circle (.05cm);
% ok
\draw[fill=black] (2cm,0cm) circle (.2cm);%%
% ok
\draw[fill=black] (3cm,1.73cm) circle (.2cm);%%
\draw[fill=black] (4cm,0cm) circle (.05cm);%%
\draw[fill=black] (6cm,0cm) circle (.05cm);%%
\draw[fill=black] (8cm,0cm) circle (.05cm);
\draw[fill=black] (-1cm,1.73cm) circle (.05cm);
% ok
\draw[fill=black] (1cm,1.73cm) circle (.2cm);%%
\draw[fill=black] (3cm,1.73cm) circle (.05cm);
\draw[fill=black] (5cm,1.73cm) circle (.05cm);
\draw[fill=black] (7cm,1.73cm) circle (.05cm);%%
\draw[fill=black] (0cm,3.46cm) circle (.05cm);
\draw[fill=black] (2cm,3.46cm) circle (.05cm);%%
\draw[fill=black] (4cm,3.46cm) circle (.05cm);
\draw[fill=black] (6cm,3.46cm) circle (.05cm);%%
\draw[fill=black] (1cm,5.20cm) circle (.05cm);
\draw[fill=black] (3cm,5.20cm) circle (.05cm);%%
\draw[fill=black] (5cm,5.20cm) circle (.05cm);%%
\draw[fill=black] (2cm,6.92cm) circle (.05cm);
\draw[fill=black] (4cm,6.93cm) circle (.05cm);
\draw[dashed] (30: 0 mm) -- (60: 90 mm);
\draw[dashed] (0: 0 mm) -- (0: 90 mm);
\draw[dashed] (0: 40 mm) -- (60: 40 mm);
\def\x{0.98};
\draw[-{Latex[length=3mm,width=2mm]}] (0,0) -- (3cm*\x,-1.73cm*\x);
\draw[-{Latex[length=3mm,width=2mm]}] (0mm,0mm) -- (0mm,3.46cm*\x);
%%
%% LABELS
%%
\node[] at (-0.8, 3.8) {$\alpha_1$};
\node[] at (3.8, -1.5) {$\alpha_2$};
\node[] at (0.9, 2.6) {$\varpi_1$};
\node[] at (2.5, -0.7) {$\varpi_2$};
\end{tikzpicture}
  \caption{The set $P_+^{q, \text{subreg}}$ of ``subregular weights'' of level $q=2$ for $\g = \mathfrak{sl}_3$}
  \label{fig:subreg}
\end{figure}
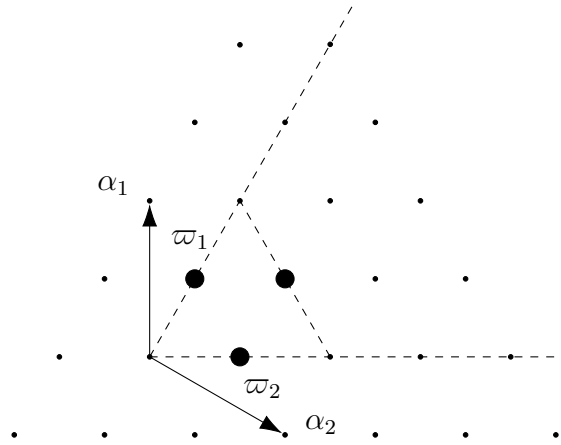

The characters $\chi_{L(\la)}(\tau, x)$ have the same expression as before. In this case $\n_+ = \n \oplus \g_{\al_*}$, so that the cancellation between the character of $\Fch(\n)$ and the Weyl denominator is less perfect than in the principal case. Indeed we obtain
\begin{align}\label{eq:indeterminate}
\chi_{C(L(\lambda))}(\tau, x) = \frac{1}{\prod_{n=1}^{\infty} (1-q^n)^\ell} \cdot \frac{B_\lambda(\tau, x)}{\Theta(\tau, \alpha_*(x))}.
\end{align}
At the end of the day we have to take the limit $x \rightarrow 0$ in \eqref{eq:indeterminate}, but here we have an indeterminate expression: The theta function has a zero along the hyperplane $\alpha_*(x) = 0$, while the numerator has a zero along the hyperplane $\gamma(x) = 0$, where $\gamma \in \D_+$ is the unique root relative to which $\lambda$ is integrable. Of course when $\lambda$ is conservative these zeros cancel out and the limit is well-defined.

Let us fix a set
\[
\{\la_1, \la_2, \ldots, \la_N\}
\]
of representatives of the classes of conservative weights, modulo $\sim$. The characters of the modules $M_i = H^0_{-}(L(\lambda_i))$ are given by the limit of \eqref{eq:indeterminate}. We would like to determine the $S$-matrix
\[
\chi_{M_i}(-1/\tau) = \sum_{i'} S_{i, i'} \chi_{M_{i'}}(\tau).
\]
In principle this can be done by performing the $S$-transformation directly upon \eqref{eq:indeterminate} with $\la = \la_i$, and carefully combining terms $B_{\mu}(\tau, x)$, for $\mu \notin \{\la_1, \ldots, \la_N\}$ to obtain cancellations and a well-defined limit. In fact we may bypass much of the difficulty of finding these explicit cancellations by keeping all the indeterminate terms, and simply computing the $x \rightarrow 0$ limit using l'Hopital's rule. 

The basic approach is to fix $x \in \h$ appropriately, and take a limit along a ray $t x$ as $t \rightarrow 0$. Terms of the form
\[
\frac{1-e^{-2\pi i \gamma(tx)}}{1-e^{-2\pi i \alpha_*(tx)}} \rightarrow \frac{\left<\gamma, x\right>}{\left<\alpha_*, x\right>}
\]
appear in the sum \eqref{eq:prin.S.mat}, whose derivation is otherwise similar. The upshot is that the $S$-matrix is
\begin{align}\label{eq:subreg.S.mat}
S_{(\nu, \eta), (\nu', \eta')} = \sum_{y(\alpha_*) \in \D_+} \varepsilon(y) \frac{\left<y(\alpha_*), x\right>}{\left<\alpha_*, x\right>} e^{-2\pi i \frac{p}{q} (y(\eta), \eta')} \sum_{w \in W} \varepsilon(w) e^{-2\pi i \frac{q}{p} (w(\nu), \nu')}.
\end{align}
The sum
\begin{align}\label{eq:degen.K}
\sum_{y(\alpha_*) \in \D_+} \varepsilon(y) \frac{\left<y(\alpha_*), x\right>}{\left<\alpha_*, x\right>} e^{-2\pi i \frac{p}{q} (y(\eta), \eta')}
\end{align}
appears here as a sort of degenerate analogue of the affine $S$-matrix \eqref{eq:Kac.sum}.

If we choose $p \geq h^\vee$ as small as possible (while respecting that $p$ and our fixed $q$ must be coprime), it often happens that the second sum in \eqref{eq:subreg.S.mat} reduces to a single term, and so the matrix \eqref{eq:degen.K}, like \eqref{eq:Kac.sum}, is realised as the $S$-matrix of a rational vertex algebra.

In many cases these rational vertex algebras are easily identified. For example the exceptional subregular $W$-algebra associated with $D_n$ and level $k = -h^\vee + p/q$ where $p/q = (2n-1)/(2n-3)$, is isomorphic to the Virasoro minimal model $\vir_{3, n-2}$, and \eqref{eq:degen.K} recovers the well-known $S$-matrix in this case \cite{FMS}.

Other cases are more intriguing. For example $\g$ of type $E_8$ and $p/q = 30/29$ yields a rational exceptional $W$-algebra whose category of modules has $44$ simple objects. The $S$-matrix \eqref{eq:degen.K} yields very complicated fusion rules via the Verlinde formula, the largest fusion coefficient being $92$ as it happens.

\newpage

\bibliographystyle{amsalpha}

\begin{thebibliography}{10}

\bibitem{AM95}
Drazen Adamovic and Antun Milas.
\newblock Vertex operator algebras associated to modular invariant representations for {$A^{(1)}_1$}.
\newblock {\em Mathematical Research Letters}, 2:563--575, 1995.

\bibitem{Ar.rep.I}
Tomoyuki Arakawa.
\newblock Representation theory of {$W$}-algebras.
\newblock {\em Inventiones mathematicae}, 169(2):219--320, 2007.

\bibitem{Ar.rep.II}
Tomoyuki Arakawa.
\newblock Representation theory of {$W$}-algebras, {II}.
\newblock {\em Exploring new structures and natural constructions in
  mathematical physics}, 61:51--90, 2011.

\bibitem{Arakawa.2015}
Tomoyuki Arakawa.
\newblock Associated varieties of modules over {K}ac-{M}oody algebras and {$C_2$}-cofiniteness of {$W$}-algebras.
\newblock {\em International Mathematics Research Notices},
  2015(22):11605--11666, 2015.

\bibitem{A.rat.O}
Tomoyuki Arakawa.
\newblock Rationality of admissible affine vertex algebras in the category
  {$\mathcal{O}$}.
\newblock {\em Duke Mathematical Journal}, 165(1):67--93, 2016.

\bibitem{A.minicourse}
Tomoyuki Arakawa.
\newblock Introduction to {$W$}-algebras and their representation theory.
\newblock {\em Perspectives in Lie theory}, pages 179--250, 2017.

\bibitem{AVE-CMP}
Tomoyuki Arakawa and Jethro van Ekeren.
\newblock Modularity of relatively rational vertex algebras and fusion rules of
  principal affine {$W$}-algebras.
\newblock {\em Communications in Mathematical Physics}, 370(1):205--247, 2019.

\bibitem{AVE-EMS}
Tomoyuki Arakawa and Jethro van Ekeren.
\newblock Rationality and fusion rules of exceptional {$W$}-algebras.
\newblock {\em J. Eur. Math. Soc. \emph{in press}}, 2022.

\bibitem{BK}
Bojko Bakalov and Alexander~A Kirillov.
\newblock {\em Lectures on tensor categories and modular functors}, volume~21
  of {\em University Lecture Series}.
\newblock American Mathematical Soc., 2001.

\bibitem{B}
Richard~E Borcherds.
\newblock Vertex algebras, {K}ac-{M}oody algebras, and the {M}onster.
\newblock {\em Proceedings of the National Academy of Sciences of the United States of America}, 83(10):3068--3071, 1986.

\bibitem{Coll.McGov}
David~H Collingwood and William~M McGovern.
\newblock {\em Nilpotent orbits in semisimple {L}ie algebra: an introduction}.
\newblock CRC Press, 1993.

\bibitem{dBT}
Jan de~Boer and Tjark Tjin.
\newblock The relation between quantum {$W$}-algebras and {L}ie algebras.
\newblock {\em Communications in mathematical physics}, 160:317--332, 1994.

\bibitem{DSK06}
Alberto De~Sole and Victor~G Kac.
\newblock Finite vs affine {$W$}-algebras.
\newblock {\em Japanese Journal of Mathematics}, 1:137--261, 2006.

\bibitem{DS}
Vladimir~G Drinfeld and Vladimir~V Sokolov.
\newblock Lie algebras and equations of {K}orteweg-de {V}ries type.
\newblock {\em Journal of Soviet mathematics}, 30:1975--2036, 1985.

\bibitem{EGNO}
Pavel Etingof, Shlomo Gelaki, Dmitri Nikshych, and Victor Ostrik.
\newblock {\em Tensor Categories}, volume 205 of {\em Mathematical Surveys and Monographs}.
\newblock American Mathematical Soc., 2015.


\bibitem{Fasquel}
Justine Fasquel.
\newblock Rationality of the exceptional {$W$}-algebras {$W_k(\mathfrak{sp}_4, f_{\text{subreg}})$} associated with subregular nilpotent elements of {$\mathfrak{sp}_4$}.
\newblock {\em Communications in Mathematical Physics}, 390(1):33--65, 2022.

\bibitem{FF}
Boris Feigin and Edward Frenkel.
\newblock Affine {K}ac-{M}oody algebras at the critical level and
  {G}elfand-{D}ikii algebras.
\newblock {\em International Journal of Modern Physics A}, 7(supp01a):197--215,
  1992.

\bibitem{Feigin.1984.semiinf}
Boris~Lvovich Feigin.
\newblock The semi-infinite homology of {K}ac-{M}oody and {V}irasoro {L}ie
  algebras.
\newblock {\em Russian Mathematical Surveys}, 39(2):155, 1984.

\bibitem{FMS}
Philippe Francesco, Pierre Mathieu, and David S{\'e}n{\'e}chal.
\newblock {\em Conformal field theory}.
\newblock Graduate Texts in Contemporary Physics. Springer-Verlag, New York,
  Inc., 1997.

\bibitem{FBZ}
Edward Frenkel and David Ben-Zvi.
\newblock {\em Vertex algebras and algebraic curves}.
\newblock Number~88 in Mathematical surveys and monographs. American
  Mathematical Soc., 2004.

\bibitem{FKW}
Edward Frenkel, Victor Kac, and Minoru Wakimoto.
\newblock Characters and fusion rules for {$W$}-algebras via quantized
  {D}rinfeld-{S}okolov reduction.
\newblock {\em Communications in mathematical physics}, 147(2):295--328, 1992.

\bibitem{FM}
IB~Frenkel and FG~Malikov.
\newblock Annihilating ideals and tilting functors.
\newblock {\em Functional Analysis and Its Applications}, 33(2):106--115, 1999.

\bibitem{FZ92}
Igor~B Frenkel and Yongchang Zhu.
\newblock Vertex operator algebras associated to representations of affine and
  {V}irasoro algebras.
\newblock {\em Duke Mathematical Journal}, 66(1):123--168, 1992.

\bibitem{GG2002}
Wee~Liang Gan and Victor Ginzburg.
\newblock Quantization of slodowy slices.
\newblock {\em International Mathematics Research Notices}, 2002(5):243--255,
  2002.

\bibitem{H.MTC}
Yi-Zhi Huang.
\newblock Rigidity and modularity of vertex tensor categories.
\newblock {\em Communications in contemporary mathematics},
  10(supp01):871--911, 2008.

\bibitem{HL1}
Yi-Zhi Huang and James Lepowsky.
\newblock A theory of tensor products for module categories for a vertex
  operator algebra, {I}.
\newblock {\em Selecta Mathematica}, 1:699--756, 1995.

\bibitem{Humphreys}
James~E Humphreys.
\newblock {\em Representations of Semisimple {L}ie Algebras in the {BGG}
  Category {$\mathcal{O}$}}, volume~94 of {\em Graduate Studies in
  Mathematics}.
\newblock American Mathematical Soc., 2008.

\bibitem{Jantz.alg}
Jens~Carsten Jantzen.
\newblock {\em Representations of algebraic groups}, volume 107 of {\em
  Mathematical Surveys and Monographs}.
\newblock American Mathematical Soc., 2003.

\bibitem{Joseph.76}
Anthony Joseph.
\newblock A characteristic variety for the primitive spectrum of a semisimple
  {L}ie algebra.
\newblock In Jacque Carmona and Mich\`{e}le Vergne, editors, {\em
  Non-Commutative Harmonic Analysis: Actes du Colloque d'Analyse Harmonique
  Non-Commutative, {M}arseille-{L}uminy, 5 au 9 {J}uillet, 1976}, volume 587 of
  {\em Lecture Notes in Mathematics}, pages 102--118. Springer, 1976.

\bibitem{KRW}
Victor Kac, Shi-Shyr Roan, and Minoru Wakimoto.
\newblock Quantum reduction for affine superalgebras.
\newblock {\em Communications in Mathematical Physics}, 241(2-3):307--342, 2003.

\bibitem{Kac.IDLA}
Victor~G Kac.
\newblock {\em Infinite-dimensional {L}ie algebras}.
\newblock Cambridge university press, 1990.

\bibitem{Kac74}
Victor~G Kac.
\newblock Infinite-dimensional {L}ie algebras and {D}edekind's {$\eta$}-function.
\newblock {\em J. Functional Anal. Appl.}, 8:68--70, 1974.

\bibitem{Kac.VA}
Victor~G Kac.
\newblock {\em Vertex algebras for beginners}.
\newblock Number~10 in University lecture series. American Mathematical Soc.,
  1998.

\bibitem{KK}
Victor~G Kac and David~A Kazhdan.
\newblock Structure of representations with highest weight of infinite-dimensional {L}ie algebras.
\newblock {\em Advances in Mathematics}, 34(1):97--108, 1979.



\bibitem{KP84}
Victor G Kac and Dale Peterson.
\newblock Infinite-dimensional {L}ie algebras, theta functions and modular forms.
\newblock {\em Advances in Mathematics},
  53(2):125--264, 1984.

\bibitem{KRR}
Victor~G Kac, Ashok~K Raina, and Natasha Rozhkovskaya.
\newblock {\em Bombay lectures on highest weight representations of infinite
  dimensional Lie algebras}, volume~29 of {\em Advanced series in mathematical
  physics}.
\newblock World scientific, 2013.

\bibitem{KW88}
Victor~G Kac and Minoru Wakimoto.
\newblock Modular invariant representations of infinite-dimensional {L}ie
  algebras and superalgebras.
\newblock {\em Proceedings of the National Academy of Sciences},
  85(14):4956--4960, 1988.

\bibitem{KW90}
Victor~G Kac and Minoru Wakimoto.
\newblock Branching functions for winding subalgebras and tensor products.
\newblock {\em Topics in Computational Algebra}, pages 3--39, 1990.

\bibitem{KW09}
Victor~G Kac and Minoru Wakimoto.
\newblock On rationality of {$W$}-algebras.
\newblock {\em Transformation Groups}, 13(3-4):671--713, 2008.

\bibitem{KW04}
Victor~G Kac, Minoru Wakimoto, et~al.
\newblock Quantum reduction and representation theory of superconformal
  algebras.
\newblock {\em Advances in Mathematics}, 185(2):400--458, 2004.

\bibitem{KL1}
David~A Kazhdan and George Lusztig.
\newblock Tensor structures arising from affine {L}ie algebras. {I}.
\newblock {\em Journal of the American Mathematical Society}, 6(4):905--947,
  1993.

\bibitem{KL2}
David~A Kazhdan and George Lusztig.
\newblock Tensor structures arising from affine {L}ie algebras. {II}.
\newblock {\em Journal of the American Mathematical Society}, 6(4):949--1011,
  1993.

\bibitem{KL3}
David~A Kazhdan and George Lusztig.
\newblock Tensor structures arising from affine {L}ie algebras. {III}.
\newblock {\em Journal of the American Mathematical Society}, 6:335--381, 1993.

\bibitem{KL4}
David~A Kazhdan and George Lusztig.
\newblock Tensor structures arising from affine {L}ie algebras. {IV}.
\newblock {\em Journal of the American Mathematical Society}, 7(2):383--453,
  1994.

\bibitem{KS}
Bertram Kostant and Shlomo Sternberg.
\newblock Symplectic reduction, {BRS} cohomology, and infinite-dimensional
  {C}lifford algebras.
\newblock {\em Annals of Physics}, 176(1):49--113, 1987.

\bibitem{Li.local.sys}
Hai-Sheng Li.
\newblock Local systems of vertex operators, vertex superalgebras and modules.
\newblock {\em Journal of Pure and Applied Algebra}, 109(2):143--195, 1996.

\bibitem{Losev}
Ivan Losev.
\newblock Finite dimensional representations of {$W$}-algebras.
\newblock {\em Duke Mathematical Journal}, 159(1), 2011.

\bibitem{Moody.Pianzola}
Robert~V Moody and Arturo Pianzola.
\newblock {\em Lie algebras with triangular decompositions}, volume~26 of {\em
  Canadian Mathematical Society Series of monographs and advanced texts}.
\newblock Wiley New York, 1995.



\bibitem{McRae-preprint}
Robert McRae.
\newblock On rationality for {$C_2$}-cofinite vertex operator algebras.
\newblock Preprint. {\em Retrived from arXiv:2108.01898}, 2021.


\bibitem{MS}
Gregory Moore and Nathan Seiberg.
\newblock Classical and quantum conformal field theory.
\newblock {\em Communications in Mathematical Physics}, 123:177-254, 1989.

\bibitem{Verlinde}
Erik Verlinde.
\newblock Fusion rules and modular transformations in {2D} conformal field theory.
\newblock {\em Nuclear Physics B}, 300:360-376, 1988.


\bibitem{zhu96}
Yongchang Zhu.
\newblock Modular invariance of characters of vertex operator algebras.
\newblock {\em Journal of the American Mathematical Society}, 9(1):237--302,
  1996.

\end{thebibliography}

\end{document}